\documentclass[pdflatex,sn-mathphys-num]{sn-jnl} 

\usepackage{graphicx}
\usepackage{multirow}
\usepackage{amsmath,amssymb,amsfonts}
\usepackage{amsthm}
\usepackage{mathrsfs}
\usepackage[title]{appendix}
\usepackage{xcolor}
\usepackage{textcomp}
\usepackage{manyfoot}
\usepackage{booktabs}
\usepackage{algorithm}
\usepackage{algorithmicx}
\usepackage{algpseudocode}
\usepackage{listings}

\usepackage{physics}

\theoremstyle{thmstyleone}

\newtheorem{property}{Property}
\newtheorem{procedure}{Procedure}
\theoremstyle{thmstyletwo}
\newtheorem{example}{Example}
\newtheorem{remark}{Remark}

\theoremstyle{thmstylethree}

\begin{document}

\title[Article Title]{A Sweeping Positivity-Preserving High Order Finite Difference WENO Scheme for Euler Equations}

\author[1]{\fnm{D. Chloe} \sur{Griffin}}\email{danielle\_griffin@brown.edu}

\author*[1]{\fnm{Chi-Wang} \sur{Shu}}\email{chi-wang\_shu@brown.edu}

\affil[1]{\orgdiv{Division of Applied Mathematics}, \orgname{Brown University}, \city{Providence}, \postcode{RI 02906}, \country{USA}}

\abstract{We develop a simple, high-order, conservative and robust positivity-preserving sweeping procedure 
for the density and the nonlinear pressure function in the compressible Euler equations. Using the scaling 
limiter in Zhang and Shu (2010) \cite{Zhang2010OnLaws}, we obtain a nontrivial extension of the scalar 
sweeping technique in Liu, Cheng, and Shu (2016) \cite{Liu2017AApproximations} for the positivity of pressure. 
The sweeping procedure developed in this paper is a post-processing technique, which can be applied to
any concave functions of the conserved variables in hyperbolic conservation law systems. Thus, it has applications
beyond the Euler equations. This procedure preserves positivity and conservation of physical quantities without 
destroying the accuracy of the underlying scheme. The algorithm works for general schemes including
finite difference, finite-volume and discontinuous Galerkin methods; however, in this paper we focus on
finite-difference weighted essentially non-oscillatory (WENO) methods. We provide numerical tests of the 
fifth order finite difference WENO scheme to 
demonstrate the accuracy and robustness of the technique.}


\maketitle

\section{Introduction}

Positivity preservation for physically meaningful quantities, particularly density and pressure in 
Euler equations, is essential for computing numerical solutions to hyperbolic conservation laws. 
Failure to maintain positivity leads to non-physical distortions, ill-posedness, and numerical 
blow-ups. Similar concerns arise in other systems, where positivity must be preserved for the density in Vlasov–Boltzmann transport equations or water height in shallow water equations. Consequently, positivity-preserving methods are crucial in diverse applications, including computational fluid dynamics (CFD), plasma simulations, computational astrophysics, traffic flow, and population models \cite{Zhang2011Maximum-principle-satisfyingDevelopments}. Many numerical methods used to solve hyperbolic conservation laws do not inherently satisfy the positivity requirements of key quantities.

Finite volume (FV) and finite difference (FD) essentially non-oscillatory (ENO) schemes \cite{Harten1997UniformlyIII,SO1}
and weighted ENO (WENO) schemes \cite{Liu1994WeightedSchemes, Jiang1996EfficientSchemes, Shu2009HighProblems}, and discontinuous 
Galerkin (DG) methods \cite{CSrev}, are among the most widely used approaches for solving hyperbolic 
conservation laws and are valued for their high-order accuracy and robust shock-capturing capabilities. Finite-difference 
WENO schemes are appealing because they are simple to implement and computationally efficient, particularly for 
multidimensional problems. Their main limitation is that their application is restricted to smooth, structured, 
curvilinear meshes, which limits their application to more complex domains. Although finite-volume WENO and DG schemes 
are more costly, they provide greater flexibility because they can be applied to arbitrary triangulations and 
non-smooth meshes.  We refer to \cite{Shu2003High-orderCFD} for a thorough description and comparison of each of these 
methods. In some cases, particularly for multidimensional systems, finite difference methods are preferred,
whenever applicable, because of 
their low computational and memory costs when compared to other methods. In comparison with the finite-volume WENO 
and DG schemes, rigorously high order positivity-preserving techniques for finite-difference schemes are less mature. 

One of the most successful strategies for ensuring the positivity of high-order schemes is the scaling limiter 
introduced in \cite{Zhang2010OnLaws}.  The idea in \cite{Zhang2010OnLaws} is to assume that the cell averages
of the high order numerical solution is positive (which has to be proved separately), then the application of
the scaling limiter can bring the whole polynomial in each cell to be positive at designated quadrature points,
without losing the high order accuracy.  This limiter has been applied to many PDEs in both finite volume 
WENO \cite{Zhang2012Maximum-principle-satisfyingEquations, Cai2016Positivity-PreservingEquations} and DG 
\cite{Zhang2010OnMeshes, Zhang2012Maximum-principle-satisfyingMeshes, Zhang2013Maximum-principle-satisfyingMeshes, 
Zhang2011Positivity-preservingTerms, Xing2010Positivity-preservingEquations, Xing2013Positivity-preservingMeshes, 
Zhang2017OnEquations} setting. For these methods, the scaling limiter is used in the reconstruction and DG 
polynomials. It is extended in \cite{Zhang2012Positivity-preservingEquations} for positivity preservation of 
density and pressure in the Euler equations for high order finite difference WENO schemes. 
Simplified implementations of the scaling limiter are provided 
in \cite{Zhang2011Maximum-principle-satisfyingDevelopments} and \cite{Wang2012RobustDetonations}. The scaling 
limiter provably maintains positivity and accuracy of the underlying scheme but requires strict theoretical 
CFL restrictions in order to prove the positivity of the cell averages before the limiter can be applied. 
This framework can be used on finite difference schemes as well \cite{Zhang2012Positivity-preservingEquations}, 
but it is more restrictive and less
popular than for finite volume and DG schemes.

Another common approach to enforce positivity preserving
for finite difference methods is via flux correction. Parameterized flux limiters use a 
combination of high- and low-order fluxes, similar to those in classic flux-corrected transport (FCT) methods 
\cite{Boris1973Flux-correctedWorks, Book1975Flux-correctedMethod, Boris1976Flux-correctedAlgorithms}. 
Examples of these flux limiters include those developed in \cite{Liang2014ParametrizedLaws, Xiong2013AFlows} 
which are generalized and applied to compressible Euler equations in \cite{Xiong2016ParametrizedEquations, 
Christlieb2015HighMeshes}. Notably, \cite{Liang2014ParametrizedLaws} has provable accuracy for certain 
schemes with some restrictions on the CFL. There are other flux limiters with direct applications to 
Euler equations in a finite-difference setting, including \cite{Hu2013Positivity-preservingEquations,
Li2025Spatial-TemporalConditions,Seal2016AnEquations}. These limiters provide positivity control for both 
density and pressure. Some methods, like those in \cite{Li2018AEquations,
Ren2024PositivityModel} and \cite{Liu2017AApproximations}, provide simplified and 
inexpensive scalar modifications for the positivity preservation of the water height in shallow water 
equations and density in Euler equations, respectively. 

Our method seeks to combine the accuracy of the scaling limiter in \cite{Zhang2010OnLaws} with the 
simplicity of \cite{Liu2017AApproximations} as an alternative positivity-preserving procedure for density, and especially for nonlinear pressure, in the Euler systems.
The sweeping procedure in \cite{Liu2017AApproximations} provides a bound-preserving technique for 
conserved scalar quantities that is independent of the underlying scheme. It is used as a 
post-processing step, with examples including one- and two-dimensional high-order 
finite-difference and finite-volume methods and a Fourier spectral method. The method guarantees 
the enforcement of bounds after a single forward and backward sweep. The procedure maintains 
the accuracy of the underlying schemes in numerical tests. In \cite{Laiu2019PositivityEquations},  
Laiu and Hauck compare several nodal and modal limiters, including the sweeping limiter. They prove 
the optimal discrete $L^1$ convergence rate of the sweeping limiter for the spectral approximation 
of kinetic transport equations in certain geometries, as well as optimal numerical convergence 
results. 

For the Euler systems, even though
the sweeping procedure in \cite{Liu2017AApproximations} can be directly applied to density, which
is one of the conserved variables, there is a significant challenge to extend it for positivity-preserving
of pressure, or internal energy, or any nonlinear concave function of the conserved variables.
In this paper, we provide a simple, conservative, positivity-preserving sweeping procedure for 
any concave functions of the conserved variables in hyperbolic conservation laws. The limiter is used as 
a post-processing step and can be applied to many functions in general schemes including 
finite difference, finite-volume and discontinuous Galerkin methods, however we will focus on the fifth order
WENO finite difference scheme in this paper.  We emphasize again that the extension of the sweeping 
limiter in \cite{Liu2017AApproximations} to the nonlinear pressure function in the Euler equations (or to
any nonlinear concave function of the conserved variables) is nontrivial. The scaling limiter in 
\cite{Zhang2010OnLaws} provides a robust adjustment mechanism for negative points; however, it relies on 
the existence of nearby points with positive pressure. 
We provide strategies for the adjustments if the local points have negative pressure during the sweeping. 
Unlike for the scalar case, we can no longer prove that the sweeping procedure will terminate after a 
single forward and backward sweep. Nevertheless, we observe numerically that, even for highly demanding test 
cases, only one or two forward and backward sweeps are required at any given Runge-Kutta step. We apply 
the pressure sweeping technique for many of the same test cases as in \cite{Zhang2012Positivity-preservingEquations}. 
The procedure maintains fifth order accuracy and also maintains 
positivity and prevents blow-ups for robust cases, including those with high CFL numbers at the instability threshold. 

The remainder of this paper is organized as follows. In Section \ref{section: background}, we provide the 
necessary description for the Euler equations, finite-difference WENO, scalar sweeping procedure (which can
be directly applied to density), and scaling limiter. We also describe the intuition behind the sweeping 
extension. In Section \ref{section: Sweeping_Procedure}, we outline the algorithm in detail, including relevant properties of the algorithm and multidimensional extension. In Section 
\ref{section: numerical_results}, we provide numerical results for the fifth order finite difference WENO scheme
with initial conditions similar to those in \cite{Zhang2012Positivity-preservingEquations}. In Section 
\ref{section: concluding_remarks}, we provide final remarks and discuss future work. 

\section{Preliminaries}
\label{section: background}

\subsection{Euler Equations}
In this study, we propose a positivity-preserving sweeping technique for the Euler equations. The one-dimensional version for the perfect gas is given by 
\begin{equation}
    \mathbf{u}_t + \textbf{f}(\mathbf{u})_x = 0,  \quad t \geq 0,\quad x \in \mathbb{R},
    \label{eq: 1D_Euler}
\end{equation} 

\[\mathbf{u} = \begin{pmatrix}
    \rho \\
    m \\
    E
\end{pmatrix}, \quad \bf{f}\bf{(u)} = \begin{pmatrix}m\\ 
\rho v^2 + p,\\
(E+p)v\end{pmatrix}\]

and
\[m = \rho v,  \quad p = (\gamma - 1)\rho e,\quad E = \frac12 \rho v^2 + \rho e,\]
where $\rho$ is density, $p$ is the pressure, $m$ is the momentum, $v$ is the velocity, $E$ is total energy, and $e$ is the internal energy. Here, $\gamma > 1$ is a constant, where $\gamma = 1.4$ is chosen for air. We consider the pressure function $p(\textbf{u}) = (\gamma - 1)(E - \frac12 \frac{m^2}{\rho})$. If $\rho > 0,$ $p(\textbf{u})$ is a concave function due to its negative semi-definite Hessian matrix. Let us consider a sequence of discrete points, $\{\textbf{u}_j\}_{j = 1}^n$. We define the 
average conserved quantity $\bar{\textbf{u}} = \sum_{j = 1}^n \frac{1}{n} \textbf{u}_j.$ Let us consider two arbitrary 
points $\textbf{u}_1$ and $\textbf{u}_2.$ From Jensen's inequality, we know that,
for $0 \le t \le 1$, 
\begin{equation} 
\label{eq: Jensen's}
p(t\textbf{u}_1 + (1 - t)\textbf{u}_2) \geq tp(\textbf{u}_1) + (1 - t)p(\textbf{u}_2), \quad \text{ if } \rho_1, \rho_2 > 0. \end{equation} The pressure sweeping algorithm is designed to ensure that all discrete point values have positive pressure. For robustness, we use an epsilon threshold for positivity, where 
$\epsilon \ll 1$. We use $\epsilon = 10^{-13}$ in our experiments. Following the notation in \cite{Zhang2012Positivity-preservingEquations}, we aim to ensure that all points are in the admissible set 

\begin{equation}
\label{eq: admissible_set}
\mathcal{G} = \left\{\bf{u} = \begin{pmatrix}
    \rho\\m\\E
\end{pmatrix} \middle\vert \rho > \epsilon \quad \text{and} \quad p = (\gamma - 1)(E - \frac12 \frac{m^2}{\rho}) >\epsilon \right\}.
\end{equation}
We note that $\mathcal{G}$ is convex. For initial value problems, we have a natural assumption that $p(\textbf{u}^0_j)>\epsilon$ for $j = 1,...,n$, where the superscript denotes the time level. From Jensen's inequality, this implies that for a 
conservative scheme and periodic or compact boundary condition, we have 
\begin{equation}
\label{eq: pressure_assumption}
p(\bar{\textbf{u}}^k) = p\left(\sum_{j = 1}^n \frac{1}{n} \textbf{u}^k_j\right) =  p\left(\sum_{j = 1}^n \frac{1}{n} \textbf{u}^0_j\right) 
\geq \sum_{j = 1}^n \frac{1}{n}p(\textbf{u}^0_j)>\epsilon\end{equation} 
at arbitrary time level $k$. Thus, it is a natural 
assumption that $p(\bar{\textbf{u}})$ is positive for all time steps. Without positivity treatment, we 
cannot assume that $\overline{p(\textbf{u})} = \sum_{j = 1}^n \frac{1}{n}p(\textbf{u}_j)>\epsilon$ once time 
stepping begins. 

Positivity of density and pressure is vital in maintaining the stability of numerical discretizations of the conservation law in (\ref{eq: 1D_Euler}). Without positivity, the linearized system becomes ill-posed, and blow-ups may 
occur in the computation. Positivity is particularly important in preventing the computation of the sound speed, $c = \sqrt{\gamma p/\rho}$ from producing imaginary values. In this study, we focus on the fifth order finite difference 
WENO as our discretization method, although our algorithm can be applied to finite-volume, finite-difference, and DG methods. 

\subsection{Finite difference WENO}

Here, we summarize the fifth order finite difference WENO scheme for (\ref{eq: 1D_Euler}) used in our applications. 
We summarize the scheme as given in \cite{Shu1998EssentiallyLaws}. We use the same implementation as  \cite{Zhang2012Positivity-preservingEquations} with a different positivity treatment. We first define a uniform one-dimensional mesh with nodes $x_i$ and mesh size $\Delta x = x_{i + 1} - x_i, $ and half-points $ x_{i + \frac12} =  \frac12(x_{i + 1} - x_i)$. We write $\textbf{u}_i(t)$ as the numerical approximation of the exact solution $\textbf{u}(x_i,t)$. 
To discretize (\ref{eq: 1D_Euler}), we use a conservative approximation 
\begin{equation}
\label{eq: conservative approx}\frac{d\textbf{u}_i(t)}{dt} = L(\textbf{u})_i = -\frac{1}{\Delta x}\left(\hat{\textbf{f}}_{i + \frac12} - \hat{\textbf{f}}_{i - \frac12}\right).\end{equation}
For the case of Euler equations, we use Lax-Friedrichs (Rusanov) splitting to ensure upwind bias,
\[
\mathbf{f}^{\pm}(\mathbf{u}) = \frac{1}{2} \left( \mathbf{u} \pm \frac{\mathbf{f}(\mathbf{u})}{\alpha} \right).\]
We compute $\alpha = \max(|\textbf{v}| + c)$ where $\textbf{v}$ is the velocity and $c$ is the sound speed. We take the maximum globally as a constant in each Runge-Kutta (RK) stage. Alternatively, the $\alpha$ computation can be performed locally in the WENO reconstruction stencil. 
At each half-point $x_{i + \frac12},$ we
\begin{enumerate}
    \item Apply flux splitting: 
    \[\mathbf{f}_i^{\pm} = \frac12\left( \mathbf{u}_i \pm \frac{\mathbf{f}( \mathbf{u}_i)}{\alpha}\right)\]
    \item Compute the Roe average matrix:
Let $\textbf{A}_{i + \frac12}$ be the Roe average matrix \cite{Roe1981ApproximateSchemes}. Compute the left and right eigenvector matrices $\textbf{L}_{i + \frac12}$ and $\textbf{R}_{i + \frac12}$ such that $\textbf{A} = \textbf{R}\mathbf{\Lambda} \textbf{L}$ where $\mathbf{\Lambda}$ is the diagonal matrix such that diagonal entries are eigenvalues of \textbf{A}.
    \item Project into characteristic fields: For each $j$ in a neighborhood of $i$, find 
    \[\mathbf{F}_j^{\pm} = \textbf{L}_{i + \frac12}\mathbf{f}^\pm_j.\]
    \item Reconstruct with fifth order WENO: Apply fifth order WENO reconstruction to each component of $\mathbf{F}^{\pm}_j$ to obtain nodal values
    $\hat{\mathbf{F}}^{\pm}_{i + \frac12}$,
 \[  \hat{\mathbf{\mathbf{F}}}^{\pm}_{i+1/2} = \text{WENO}\big(\mathbf{F}^{\pm}_{i-r+1}, \dots, \mathbf{F}^{\pm}_{i+r}\big),\]
 where $r = 3$ for the fifth-order WENO reconstruction.

    \item Project back into physical space:\\
   \[ \hat{\mathbf{f}}^{\pm}_{i + \frac12}  = \textbf{R}_{i + \frac12} \hat{\mathbf{F}}^{\pm}_{i + \frac12}.\]
   
    \item Form fluxes and plug into conservative scheme:
    \[\hat{\mathbf{f}}_{i + \frac12} = \alpha(\hat{\mathbf{f}}_{i + \frac12}^+ - \hat{\mathbf{f}}_{i + \frac12}^-)\]
and obtain $L(\mathbf{u}_i)$ as in \eqref{eq: conservative approx}.
\end{enumerate}
For the fully discrete method, we use third-order Runge-Kutta time-stepping: 
\begin{equation}
\begin{aligned}
\textbf{u}^{(1)} &= \textbf{u}^n + \Delta t \, L(\textbf{u}^n), \\
\textbf{u}^{(2)} &= \frac{3}{4} \textbf{u}^n + \frac{1}{4}( \textbf{u}^{(1)} + \Delta t \, L(\textbf{u}^{(1)})), \\
\textbf{u}^{n+1} &= \frac{1}{3} \textbf{u}^n + \frac{2}{3}( \textbf{u}^{(2)} + \Delta t \, L(\textbf{u}^{(2)})).
\end{aligned}
\end{equation}
 We apply the sweeping limiter as a post-processing step after each RK stage. The limiter is only 
applied whenever there is negative pressure or density.

\subsection{Density sweeping procedure}

We use the simple density sweeping procedure in \cite{Liu2017AApproximations} and the scaling limiter in \cite{Zhang2010OnLaws} as the primary building blocks of the pressure-sweeping algorithm. In \cite{Liu2017AApproximations}, a sweeping algorithm is provided for a generic, conserved, scalar quantity in either a finite difference or finite volume implementation. To modify density in the Euler equations in a finite-difference WENO setting, we consider a collection of point values 
$\{\rho_j\}_{j = 1}^n,$ with the average density $\bar{\rho} = \sum_{j = 1}^n\frac{1}{n} \rho_ j $. Here, we rewrite the density-sweeping algorithm for clarity.

\begin{algorithm}[H]
\label{alg: density_sweep}
\caption{Positivity-preserving density sweep (Algorithm 2.1  of \cite{Liu2017AApproximations})}
\label{alg:bp-sweep}
\begin{algorithmic}
    \State $\epsilon \gets 10^{-13}$
    \For{$j=1$ : $n-1$} \Comment{Forward Sweep}
  \If{$\rho_j < \epsilon$}
    \State $\displaystyle d \gets \rho_j - \epsilon$ 
    \State $\rho_{j+1} \gets \rho_{j+1} + d$
    \State $\rho_j \gets \epsilon$
  \EndIf
\EndFor

\For{$j=n$ : $2$} \Comment{Backward Sweep}
  \If{$\rho_j < \epsilon$}
    \State $\displaystyle d \gets \rho_j - \epsilon$ 
    \State $\rho_{j-1} \gets \rho_{j-1} + d$
    \State $\rho_j \gets \epsilon$
  \EndIf
\EndFor
\end{algorithmic}
\end{algorithm}

The sweeping algorithm has several advantages. It is simple to implement and can be used as a post-processing step in many numerical methods where positivity is required. For example, the sweeping technique is used in \cite{Tsybulnik2023EfficientEquations} for the numerical solutions of the 3D-Vlasov-Maxwell system. As a global serial process, a disadvantage is that it cannot be parallelized. Nevertheless, only a single forward and backward sweep is required to preserve positivity. Thus, this approach is relatively inexpensive. In our application of the density sweeping procedure, and our extension for pressure, we use parallel computation for spatial discretization and only apply the serial function as a post-processing step in the Runge-Kutta time discretization. In the numerical tests provided in \cite{Liu2017AApproximations}, the sweeping technique preserves 
the accuracy of the underlying scheme. Recent theoretical convergence results have also been obtained. 
In \cite{Laiu2019PositivityEquations}, Laiu and Hauck evaluate the sweeping limiter's use in the spectral approximation of kinetic transport equations. Specifically, they apply sweeping to the values of the spectral approximation at certain quadrature points. They compare several nodal and modal limiters, finding that sweeping outperformed both modal limiters, as well as the clipping limiter in \cite{Shashkov2004TheLaws} in terms of accuracy. They introduce a new hierarchical version of the clipping limiter with similar accuracy. Although hierarchical clipping has the benefit of parallel implementation, it suffers from more numerical artifacts in certain test cases. Laiu and Hauck provide theoretical and numerical convergence results for the limiters. They provide optimal discrete $L^1$ convergence rates for the sweeping limiter on slab geometries, which are sharp in numerical applications. They also provide suboptimal discrete $L^1$ consistency estimates on $\mathbb{S}^2,$ although numerical results indicate that the results are not sharp. 

We choose to extend the sweeping limiter owing to its accuracy and simple implementation. However, the extension to the
nonlinear pressure function is not straightforward. The density sweeping algorithm relies on the fact that $\bar{\rho} > \epsilon$ is conserved. Thus, upon gathering all negative density in the forward sweep and redistributing the negative density in the backward sweep, there must be nothing left over at the end. Otherwise, this would violate $\bar{\rho} > \epsilon$. This is the non-technical description of the proof given in \cite{Liu2017AApproximations} that only one forward and backward sweep is
needed. We have no guarantee that the average pressure, $\overline{p(\textbf{u})} = \sum_{i = 1}^np(\textbf{u}_j) >\epsilon $. Combining the results underlying assumption \eqref{eq: pressure_assumption} and Jensen's inequality \eqref{eq: Jensen's}, concavity of pressure and conservation only ensures that 
\[\overline{p(\textbf{u})} \leq p(\overline{\textbf{u}}) \text{ and } p(\overline{\textbf{u}}) > \epsilon.\]
It is possible for $\overline{\textbf{u}} \in \mathcal{G}$ while $\textbf{u}_j \not \in \mathcal{G}$ for $j = 1,..,n$. Thus, even if we can easily project any negative value $\textbf{u}_j$ to the admissible set $\mathcal{G}$ and make conservative exchanges as in the original sweeping procedure, we cannot guarantee that only one forward and backward sweep is 
needed. Furthermore, projecting onto a convex set is expensive. Thus, we provide an alternative adjustment style using the scaling limiter in \cite{Zhang2010OnLaws}.

\subsection{Scaling limiter}

Zhang and Shu introduce the scaling limiter in \cite{Zhang2010OnLaws} as a bound preserving technique where 
reconstruction or DG polynomials are scaled linearly around the cell average (which is proved to be
positive separately) while maintaining uniform high 
order accuracy under certain CFL conditions. This strategy is an extension of the approach in 
\cite{Perthame1996OnEquations} and is applied in a finite difference setting for Euler equations in 
\cite{Zhang2012Positivity-preservingEquations}. We consider the scaling limiter for the treatment of 
pressure in the Euler equations. Assuming $\rho_1, \rho_2 > \epsilon$, let us consider $\textbf{u}_1, \textbf{u}_2 \in C$ 
such that $p(\textbf{u}_1) \leq \epsilon$ and $p(\textbf{u}_2)>\epsilon$. Thus, $\textbf{u}_2 \in \mathcal{G}$ and 
$\textbf{u}_1 \not \in \mathcal{G}$. Motivated by the sweeping limiter in \cite{Liu2017AApproximations}, we 
make a nodal adjustment of $\textbf{u}_1$ with its neighbor such that $\textbf{u}^*_1 \in \mathcal{G}$. Following the approach in \cite{Zhang2010OnMeshes}, one nodal adjustment choice is to set
\begin{equation}
\label{eq: line_move}
    \textbf{u}_1^* = t(\textbf{u}_2 - \textbf{u}_1) + \textbf{u}_1,
\end{equation}
where $t$ satisfies 
\begin{equation}
\label{eq}
    p((1 - t)\textbf{u}_1 + t\textbf{u}_2) = \epsilon.
\end{equation}Thus, $t$ provides the parameterized amount to move along the line between $\textbf{u}_1$ and $\textbf{u}_2$ for $\textbf{u}_1^*$ to be in the admissible set $\mathcal{G}$, but requires a quadratic solver. A proof is also provided in \cite{Zhang2010OnMeshes} that verifies that the limiter does not destroy the accuracy of smooth solutions in its application to discontinuous Galerkin methods. This result extends to a simplified and more robust method to find $t$ given in \cite{Wang2012RobustDetonations}. We rewrite it here in terms of our applications. Using Jensen's inequality, we solve for $t$ which satisfies 
\[\epsilon = (1-t)p(\textbf{u}_1) + tp(\textbf{u}_2) \leq p((1-t)\textbf{u}_1 + t(\textbf{u}_2)) = p(\textbf{u}_1^*).\]
Thus, for $\rho_1, \rho_2 > \epsilon$ and $p(\textbf{u}_1)\leq\epsilon,$ we ensure $p(\textbf{u}_1^*) \geq \epsilon$ by setting
\begin{equation}
\label{eq: t_pos_adjustment}
t = \frac{p(\textbf{u}_1) - \epsilon}{p(\textbf{u}_1) - p(\textbf{u}_2)}.\end{equation}
We note that the more robust definition is necessary for problems involving blast waves. Thus far, the sweeping limiter only provides an adjustment protocol for $p(\textbf{u}_1) \leq \epsilon$  whenever we also have 
$p(\textbf{u}_2) > \epsilon $. We restrict exchanges to neighboring nodes in hopes of improving accuracy. Furthermore, we are not guaranteed that there is another positive point to exchange with, even if we consider all points. Therefore, we must provide an alternative mechanism to ensure that all points are in $\mathcal{G}$. For this, we use the scaling limiter to obtain an estimate of the distance to the convex set $\mathcal{G}$ from a negative point. The only guaranteed point 
in the set is $\bar{\textbf{u}}$. Thus, we estimate that distance by solving
\begin{equation}
\label{eq: t_pos_adjustment2}
t_1 = \frac{p(\textbf{u}_1) - \epsilon}{p(\textbf{u}_1) - p(\bar{\textbf{u}})}\end{equation}
and scale along the line between $\textbf{u}_1$ and $\textbf{u}_2$. Thus, we aim to move the necessary distance to reach the set $\mathcal{G}$ if the line between $\textbf{u}_1$ and $\textbf{u}_2$ intersects with the set. If the line does not intersect, there is no way for this nodal adjustment to indicate $\textbf{u}^*_1 \in \mathcal{G}$. In this case, the best we can do is 
to ensure that the two points move closer together, thereby increasing their combined pressure. Thus, we include a minimum of $1/4$. In this way, the distance between negative points will at most be halved. If they are already close to $\mathcal{G},$ a smaller distance may be provided by the scaling limiter estimate. The full definition for $t$ in this case is 
\[ t = \min\left(t_1\left(\frac{\norm{\textbf{u}_1 - \bar{\textbf{u}}}}{\norm{\textbf{u}_1 - \textbf{u}_2}}\right), \frac{1}{4}\right).\]
We provide a full description of the algorithm in the next section.

\section{Conservative positivity-preserving sweeping procedure for pressure}
\label{section: Sweeping_Procedure}

Here, we describe the pressure sweeping algorithm with positivity-preserving properties. Given $\{\textbf{u}_j\}_{j = 1}^n$ with 
$p(\bar{\textbf{u}}) = p(\sum_{j = 1}^n \frac{1}{n} \textbf{u}_j) > \epsilon$, the pressure sweeping procedure continues until positivity of all $\rho_j \in \{\rho_j\}_{j = 1}^n$ and $p(\textbf{u}_j) \in  \{p(\textbf{u}_j)\}_{j = 1}^n$ is achieved. In other words, it provides a method to adjust all points towards admissible set $\mathcal{G}$ given in (\ref{eq: admissible_set}) under the natural assumption in (\ref{eq: pressure_assumption}). Furthermore, it maintains the conservation of 
$\textbf{u}$ such that $\bar{\textbf{u}} = \sum_{j = 1}^n \frac{1}{n} \textbf{u}_j$ is unchanged. 
To correct a ``negative point" $\textbf{u}_j$, or a point where $p(\textbf{u}_j) < \epsilon,$ we find an adjusted point $\textbf{u}_j^*$ using a convex combination of $\textbf{u}_j$ with its neighbor. As in the original sweeping algorithm, we make a nodal adjustment with the next point in the direction of sweep. We only adjust the negative points, and the nodal adjustment differs if the neighboring point is positive or negative. 

If the neighboring point is positive, $\textbf{u}_j^*$ is adjusted using the scaling limiter as in \eqref{eq: line_move} with $t$ from (\ref{eq: t_pos_adjustment}). If the neighboring point is negative, we use a similar estimate to find the distance to the convex set $\mathcal{G}$ using $\bar{\textbf{u}}$. We use this distance as an estimate of how much to move along a line between the negative point and its neighbor. We take a minimum of this scaled distance with $1/4$. In this way, the distance between negative points is at most halved. When points are close to the admissible set, the estimate by the scaling limiter may provide a better estimate of the distance needed for $\textbf{u}_j^* \in \mathcal{G}$. 

In practice, $\{\textbf{u}_j\}$ may be a sequence of vectors with point values or cell averages,
but we provide the algorithm for point values in a finite difference scheme. 
Throughout this work, we denote by $\norm{\cdot}$ the discrete $L^2$ norm in $\mathbb{R}^3$, defined for any vector $\mathbf{x} = (x_1, x_2, x_3) \in \mathbb{R}^3$ as
\[
\norm{\mathbf{x}} = \sqrt{x_1^2 + x_2^2 + x_3^2}.
\]The algorithm is as follows:

\begin{procedure}(A conservative positivity-preserving sweeping procedure)
\label{algorithm: pos_preserve_sweep}
\begin{enumerate}
    \item Apply Algorithm \ref{alg: density_sweep} to $\{\rho_j\}_{j = 1}^n$,
    \item Apply the positivity preserving sweep in Algorithm \ref{alg: positivity_sweep} using the nodal adjustment in Algorithm \ref{alg: node_adjust}.

    \begin{algorithm}[H]
    \caption{Positivity Preserving Sweep}
     \label{alg: positivity_sweep}
        \begin{algorithmic}
        \State $\epsilon \gets 10^{-13}$
        \While{$\min(p(\textbf{u}))<\epsilon$}
        \For{$j = 1 : n - 1$} \Comment{Forward Sweep}
        \If{$ p(\textbf{u}_j) < \epsilon $}
        \State ${\textbf{u}}^*_j \gets \text{NodeAdjust}(\textbf{u}_j, \textbf{u}_{j + 1}, \bar{\textbf{u}}, \epsilon)$ 
        \State $\textbf{d} \gets \textbf{u}_j - \textbf{u}^*_j$
        \State $\textbf{u}^*_{j + 1} \gets \textbf{u}_{j + 1} + \textbf{d}$
        \State $\textbf{u}_j \gets \textbf{u}^*_j$
        \State $\textbf{u}_{j + 1} \gets \textbf{u}^*_{j + 1}$
        \EndIf   
            \EndFor
        \For{$j = n : 2$} \Comment{Backward Sweep}
        \If{$ p(\textbf{u}_j) < \epsilon $}
        \State $\textbf{u}^*_j \gets \text{NodeAdjust}(\textbf{u}_j, \textbf{u}_{j - 1}, \bar{\textbf{u}},\epsilon)$
        \State $\textbf{d} \gets \textbf{u}_j - \textbf{u}^*_j$
        \State $\textbf{u}^*_{j - 1} \gets \textbf{u}_{j - 1} + \textbf{d}$
        \State $\textbf{u}_j \gets \textbf{u}^*_j$
        \State $\textbf{u}_{j - 1} \gets \textbf{u}^*_{j - 1}$
        \EndIf   
        \EndFor
        \EndWhile
        \end{algorithmic}
        \end{algorithm}
 \end{enumerate}  
The following function provides a node adjustment protocol.

    \begin{algorithm}[H]
     \caption{$\textbf{u}^*_1$ = NodeAdjust($\textbf{u}_1, \textbf{u}_2, \bar{\textbf{u}}, \epsilon)$}
      \label{alg: node_adjust}
         \begin{algorithmic}
        \If{$ p(\textbf{u}_2) > \epsilon $}
            \State $t \gets \left(\frac{p(\textbf{u}_1) - \epsilon}{p(\textbf{u}_1) - p(\textbf{u}_2)}\right)$
         \Else
            \State $t_1 \gets \frac{p(\textbf{u}_1) - \epsilon}{(p(\textbf{u}_1) - p(\bar{\textbf{u}}))}$
            \State $t \gets \min\left(t_1\left(\frac{\norm{\textbf{u}_1 - \bar{\textbf{u}}}}{\norm{\textbf{u}_1 - \textbf{u}_2}}\right), 1/4\right)$
        \EndIf   

        \State $\textbf{u}^*_1 \gets (1 - t)\textbf{u}_1 + t\textbf{u}_2$
        \end{algorithmic}
        \end{algorithm}
 \end{procedure}
We describe several useful properties of the algorithm.

\begin{property}[Conservation]
\label{Property: conservation}
    First, we note that none of the operations in Algorithm \ref{algorithm: pos_preserve_sweep} change $\bar{\textbf{u}}$. 
Denoting $\textbf{u}_n$ as a node where $p(\textbf{u}_n) < \epsilon$ and $\textbf{u}_p$ as the neighboring node, we have
    \[\textbf{u}_n^* + \textbf{u}_p^* = \textbf{u}_n^* + (\textbf{u}_p + (\textbf{u}_n - \textbf{u}_n^*)) = \textbf{u}_n + \textbf{u}_{p}. \]
This implies that $\bar{\textbf{u}}$ is not changed and hence we still have $p(\overline{\textbf{u}}) > \epsilon$. 
\end{property}

\begin{property}[Monotonic Improvement]
Any operation will only 
decrease the distance between points and increase the average of pressure, $\overline{p(\textbf{u})} = \sum_{i = 1}^np(\textbf{u}_j)$. We first rewrite
    \begin{align*}
        \textbf{u}_p^* =& \textbf{u}_p + (\textbf{u}_n - \textbf{u}_n^*)\\
              =& \textbf{u}_p + \textbf{u}_n - ((1 - t)\textbf{u}_n + t\textbf{u}_p))\\
              =& t\textbf{u}_n + (1 - t)\textbf{u}_p.    \end{align*}
    Using the definition on $\textbf{u}_p^*$ and Cauchy-Schwarz, we find a bound on the adjusted distance between points as 
    \begin{align*}
    \norm{\textbf{u}_p^* - \textbf{u}_n^*} = &\norm{\textbf{u}_p + (\textbf{u}_n - \textbf{u}_n^*) -\textbf{u}_n^*}\\
    =&\norm{\textbf{u}_p + (\textbf{u}_n - ((1 - t)\textbf{u}_n + t\textbf{u}_p)) - ((1 - t)\textbf{u}_n + t\textbf{u}_p)}\\
    =&\norm{(1 - 2t)\textbf{u}_p - (1 - 2t)\textbf{u}_n}\\
    =&\norm{(1 - 2t)(\textbf{u}_p - \textbf{u}_n)}\\
    \leq & |1 - 2t|\norm{(\textbf{u}_p - \textbf{u}_n)}.
    \end{align*}
    Letting $a = p(\textbf{u}_1) - \epsilon$ and $b = p(\textbf{u}_2) + \epsilon$ or $b = p(\bar{\textbf{u}}) + \epsilon$, it is easy to see that $t \in (0,1)$ from the simple fact that
    \[
a < 0 < b \;\implies\; 0 < \frac{a}{a - b} < 1.
\]
 Thus,
    \[|1 - 2t| \leq 1 \implies \norm{\textbf{u}_p^* - \textbf{u}_n^*} \leq \norm{\textbf{u}_p - \textbf{u}_n}.\]
Furthermore, the concavity of the pressure implies that as points move closer together, their combined pressure increases. Using Jensen's inequality, we find that 
      \begin{align*}
  p(\textbf{u}_p^*) +  p(\textbf{u}_n^*) =& p(t\textbf{u}_n + (1 - t)\textbf{u}_p) +  p((1 - t)\textbf{u}_n + t\textbf{u}_p)\\
        \geq& tp(\textbf{u}_n) + p(\textbf{u}_p) - tp(\textbf{u}_p) + p(\textbf{u}_n)- tp(\textbf{u}_n) + tp(\textbf{u}_p)\\
        =& p(\textbf{u}_p) + p(\textbf{u}_n).
    \end{align*}
It is clear that with successive application of the sweeping limiter, the average value of pressure, that is, 
$\frac 1n \sum_{i = 1}^np(\textbf{u}_j)$ is monotonically increasing and bounded above by $p(\bar{\textbf{u}})$. Thus, sweeping may only improve the conditions for pressure positivity. 
\end{property}

\begin{property}[Variance Reduction]
      As the points get closer together, the total variance decreases. As the sweeping procedure is continuously applied, all elements monotonically move toward the set $\mathcal{G}$ and expected value, $\bar{\textbf{u}}$.  To show this, we consider how the total variance $\mathbb{V} = \sum_{j = 1}^n \norm{\textbf{u}_j - \bar{\textbf{u}}}^2$ changes upon each combined nodal adjustment. We consider the local variance for $\textbf{u}_p$ and $\textbf{u}_n$ defined as 
\[\mathbb{V}_{n,p} := \norm{\mathbb{V}_n}^2 + \norm{\mathbb{V}_p}^2 := \norm{\textbf{u}_n - \bar{\textbf{u}}}^2 + \norm{\textbf{u}_p- \bar{\textbf{u}}}^2, \]
and how it changes in relation to 
\[\mathbb{V}_{n,p}^* := \norm{\mathbb{V}_n^*}^2 + \norm{\mathbb{V}_p^*}^2 := \norm{\textbf{u}_n^* - \bar{\textbf{u}}}^2 + \norm{\textbf{u}_p^*- \bar{\textbf{u}}}^2. \]
Through a simple calculation, we see that 
\[\mathbb{V}_n^* = \textbf{u}_n^* - \bar{\textbf{u}} = (1 - t)\textbf{u}_n + t\textbf{u}_p - \bar{\textbf{u}} = (1 - t)(\textbf{u}_n - \bar{\textbf{u}}) + t(\textbf{u}_p - \bar{\textbf{u}}) = (1 - t)\mathbb{V}_n + t\mathbb{V}_p .\]
Through the same calculation, we find $\mathbb{V}_p^*  = \textbf{u}_p^* - \bar{\textbf{u}} = t\mathbb{V}_n + (1 - t)\mathbb{V}_p$. We expand
\begin{align*}
    \mathbb{V}_{n,p}^* =& \norm{\mathbb{V}_n^*}^2 + \norm{\mathbb{V}_p^*}^2 \\
    =& \norm{(1 - t)\mathbb{V}_n + t\mathbb{V}_p }^2 + \norm{t\mathbb{V}_n + (1 - t)\mathbb{V}_p}^2\\
    =& (1 - t)\norm{\mathbb{V}_n}^2 + t^2\norm{\mathbb{V}_p}^2 + 2t(1 - t)\mathbb{V}_n\mathbb{V}_p + \\
 + &t^2 \norm{\mathbb{V}_n}^2 + (1 - t)^2 \norm{\mathbb{V}_p}^2 + 2t(1 - t)\mathbb{V}_n\mathbb{V}_p\\
 =& [(1 - t)^2 + t^2]\norm{\mathbb{V}_n}^2 + [t^2 + (1 - t)^2]\norm{\mathbb{V}_p}^2 + 4t(1 - t)\mathbb{V}_n\mathbb{V}_p\\
 =& (2t(t - 1) + 1)(\norm{\mathbb{V}_n}^2 + \norm{\mathbb{V}_p}^2) + 4t(1 - t)\mathbb{V}_n\mathbb{V}_p.
\end{align*}
We use the fact that 
\[\norm{\mathbb{V}_p - \mathbb{V}_n}^2 = \norm{\mathbb{V}_p}^2 + \norm{\mathbb{V}_n}^2 - 2\mathbb{V}_n\mathbb{V}_p \implies 2\mathbb{V}_n\mathbb{V}_p = \norm{\mathbb{V}_p}^2 + \norm{\mathbb{V}_n}^2 - \norm{\mathbb{V}_p - \mathbb{V}_n}.\]
Thus, 
\begin{align*}4t(1 - t)\mathbb{V}_n\mathbb{V}_p = 2t(1 - t)\left[\norm{\mathbb{V}_p}^2 + \norm{\mathbb{V}_n}^2 - \norm{\mathbb{V}_p - \mathbb{V}_n}\right].\end{align*}
Combining these results, we obtain the estimate
\begin{align*}
    \mathbb{V}_{n,p}^* =& \norm{\mathbb{V}_p}^2 + \norm{\mathbb{V}_n}^2 - 2t(1 - t)\norm{\mathbb{V}_p - \mathbb{V}_n}\\
    =& \norm{\mathbb{V}_p}^2 + \norm{\mathbb{V}_n}^2 - 2t(1 - t)\norm{\mathbf{u}_p - \mathbf{u}_n}.
\end{align*}
For $t \in (0,1),$ $0 < 2t(1 - t) < 1$. Thus, the local variance $\mathbb{V}_{p,n} < \mathbb{V}^*_{p,n}$ as long as $\mathbf{u}_p \not =  \mathbf{u}_n$. Other point values are not affected. This implies that the total variance $\mathbb{V}$ is monotonically decreasing.
\end{property} 

\begin{remark}
    While we do not provide rigorous convergence results to the set $\mathcal{G}$, we justify that pressure is increasing and variance decreases while conservation is maintained. Thus, points can only move closer to the admissible set, and negative pressure values are improved. In practice, only a few points have negative pressure, and they are likely to have positive neighbors nearby. Thus, we conjecture that the combination of the robust scaling limiter for pressure correction and the mass-redistribution technique of the sweeping limiter explains the need for only one or two forward and backward sweeps to achieve positivity in our numerical results in Section \ref{section: numerical_results}. 
\end{remark}

\subsection{Extension to multiple dimensions}

In two spatial dimensions, the Euler equations can be written in conservative form as
\begin{equation}
\frac{\partial \mathbf{w}}{\partial t}
+ \frac{\partial \textbf{f}(\mathbf{w})}{\partial x}
+ \frac{\partial \textbf{g}(\mathbf{w})}{\partial y} = 0, (x,y) \in \mathbb{R}^2, t \geq 0,
 \label{eq: 2D_Euler}
\end{equation}
where the vector of conserved variables is $\mathbf{w} = (\rho, m,n,E)^T$. Here, $\rho$ is the density, $m = \rho u$ and $n = \rho v$ are the momentum components in the $x$- and $y$-directions, respectively, and $E$ is the total energy. The flux vectors are
\begin{equation}
\textbf{f}(\mathbf{w}) =
\begin{pmatrix}
m \\[2pt]
\dfrac{m^2}{\rho} + p \\[2pt]
\dfrac{m n}{\rho} \\[2pt]
\dfrac{m}{\rho}(E+p)
\end{pmatrix}, \qquad
\textbf{g}(\mathbf{w}) =
\begin{pmatrix}
n \\[2pt]
\dfrac{m n}{\rho} \\[2pt]
\dfrac{n^2}{\rho} + p \\[2pt]
\dfrac{n}{\rho}(E+p)
\end{pmatrix}.
\end{equation}
The pressure is determined from the ideal-gas law
\[
p = (\gamma - 1)\left(E - \tfrac{1}{2}\frac{m^2+n^2}{\rho}\right),
\]
with $\gamma>1$ describing the ratio of the specific heats. As before, one can easily verify that $p$ is a concave function for $\rho > 0$. The admissible set now becomes 
\begin{equation}
\label{eq: admissible_set_2D}
\mathcal{G} = \left\{\bf{w} = \begin{pmatrix}
    \rho\\m\\n\\E
\end{pmatrix} \middle\vert \rho > \epsilon \quad \text{and} \quad 
p = (\gamma - 1)\left(E - \tfrac{1}{2}\frac{m^2+n^2}{\rho}\right) >\epsilon \right\}.
\end{equation}
which is still convex. For a two-dimensional uniform mesh with nodes $(x_i,y_j),$ we use a conservative approximation 
\[\frac{d\textbf{w}_{i,j}(t)}{dt} = L(\textbf{w})_{i,j} = -\frac{1}{\Delta x}\left(\hat{\textbf{f}}_{i + \frac12,j} - \hat{\textbf{f}}_{i - \frac12,j}\right) -\frac{1}{\Delta y}\left(\hat{\textbf{g}}_{i,j + \frac12} - \hat{\textbf{g}}_{i,j + \frac12} \right).\] 
As in the one-dimensional case, we use Lax-Friedrichs (Rusanov) splitting to ensure upwind bias,

\begin{equation}
\label{eq: conservative approx2}
\mathbf{f}^{\pm}(\mathbf{w}) = \frac{1}{2} \left( \mathbf{w} \pm \frac{\mathbf{f}(\mathbf{w})}{\alpha_1} \right)\quad  \text { and }\quad  \mathbf{g}^{\pm}(\mathbf{w}) = \frac{1}{2} \left( \mathbf{w} \pm \frac{\mathbf{g}(\mathbf{w})}{\alpha_2} \right),\end{equation}
where $\alpha_1 = \max(|u| + c)$ and $\alpha_2 = \max(|v| + c)$. We take the maximum globally as a constant in each RK stage. Alternatively, the computation can be performed locally in the WENO reconstruction stencil. 

The extension of the pressure sweeping algorithm to two-dimensions is simple. We generate indexing arrays much like in \cite{Liu2017AApproximations}. We maintain conservation over the entire computation domain. In other words, we maintain $\bar{\textbf{w}} = \sum_{i= 1}^{n_x}\sum_{j = 1}^{n_y} \frac{1}{n_xn_y}\textbf{w}_{i,j}.$ We treat $\{\textbf{w}_{k}\}_{k=1}^{n_xn_y}$ as a one-dimensional matrix. There are many ways to create this matrix while maintaining neighboring relations of the 2D array. We alternate between the two orderings given in \cite{Liu2017AApproximations}. The orderings are written as:

\begin{itemize}
    \item Sweep I: Connect mesh cells column by column: For each column, $j = 1,...,n_y$, use order $i = 1,2,...,n_x$ for odd $j,$ and order $i = n_x, n_x - 1,...,1$ for even $j.$
     \item Sweep II: Connect mesh cells row by row: For each row,  $i = 1,2,...,n_x$, use order $j = 1,...,n_y$ for odd $i,$ and order $j = n_y, n_y - 1,...,1$ for even $i.$
\end{itemize}

We alternate between Sweep I and Sweep II in our numerical tests. We count a forward and backward sweep as one full pressure sweep. Thus, at a given RK step, Sweep I is used for odd sweeps and Sweep II is used for even sweeps. This choice ensures access to neighbors in both the $x$ and $y$ directions. Nevertheless, most of the test cases in Section \ref{section: numerical_results} require a single sweep at any Runge–Kutta step. Most of our test cases need only one sweep. Thus, we expect that using only Sweep I or Sweep II would be sufficient. There is also a possibility of parallel implementation in the two-dimensional case because we allow more than a single forward and backward sweep.

\section{Numerical results}
\label{section: numerical_results}
We apply the sweeping procedure to maintain the positivity of density and pressure for both the one- and two-dimensional cases. The numerical experiments in this study are designed to mirror those in \cite{Zhang2012Positivity-preservingEquations}, demonstrating the accuracy and positivity-preserving properties of the sweeping method. We present seven numerical examples in which negative pressure occurs. In most of these cases, negative pressure causes the WENO algorithm to blow up without positivity treatment. In cases where negative values did not occur in a stable CFL range, we increased the CFL beyond the standard limits. Applying sweeping at the instability threshold, we still see key features of the numerical solution which demonstrates the robustness of the algorithm. For detailed accounts of the number of sweeps used (both in total and average/max quantities for each RK stage), see Table \ref{table: sweeping_info} and Table \ref{table: den_sweeping_info} for pressure and density respectively.

\begin{example}[1D double rarefaction]
\label{ex: 1D_rarefaction}
We first present a double rarefaction problem in one dimension. The exact solution contains a vacuum; thus, negative values may appear. We solve the 1D Euler equations given by (\ref{eq: 1D_Euler}) with the Riemann initial condition 
\[\begin{cases}
    (\rho, v, p)  = (7, -1, 0.2) & \text{if } x \leq 0,\\
    (\rho, v, p)  = (7, 1,  0.2) & \text{if } x > 0,
\end{cases}\]and $\gamma = 1.4$. The boundaries are outflow. We run our numerical test on $x = [0.5,0.5]$ with 200 points 
until the final time $T = 0.3$. We use CFL = 0.9 to produce negative pressure values. Even with a larger CFL, our results in Figure \ref{fig:1D_rarefaction} are comparable to those in \cite{Zhang2012Positivity-preservingEquations}.

\begin{figure}[h!]
    \centering
    \begin{minipage}[b]{0.45\textwidth}
        \centering
        \includegraphics[width=\textwidth]{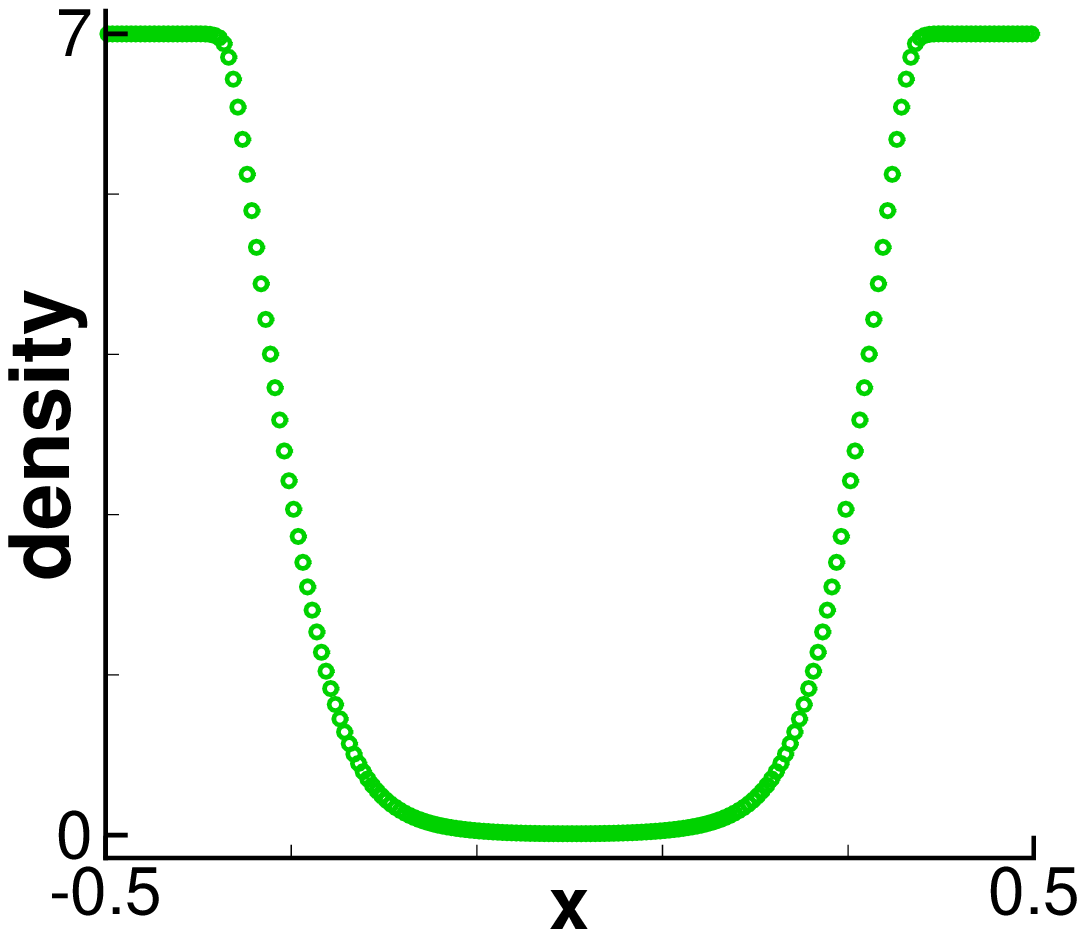}
        \\ {\small Density Plot}
    \end{minipage}
    \hfill
    \begin{minipage}[b]{0.45\textwidth}
        \centering
        \includegraphics[width=\textwidth]{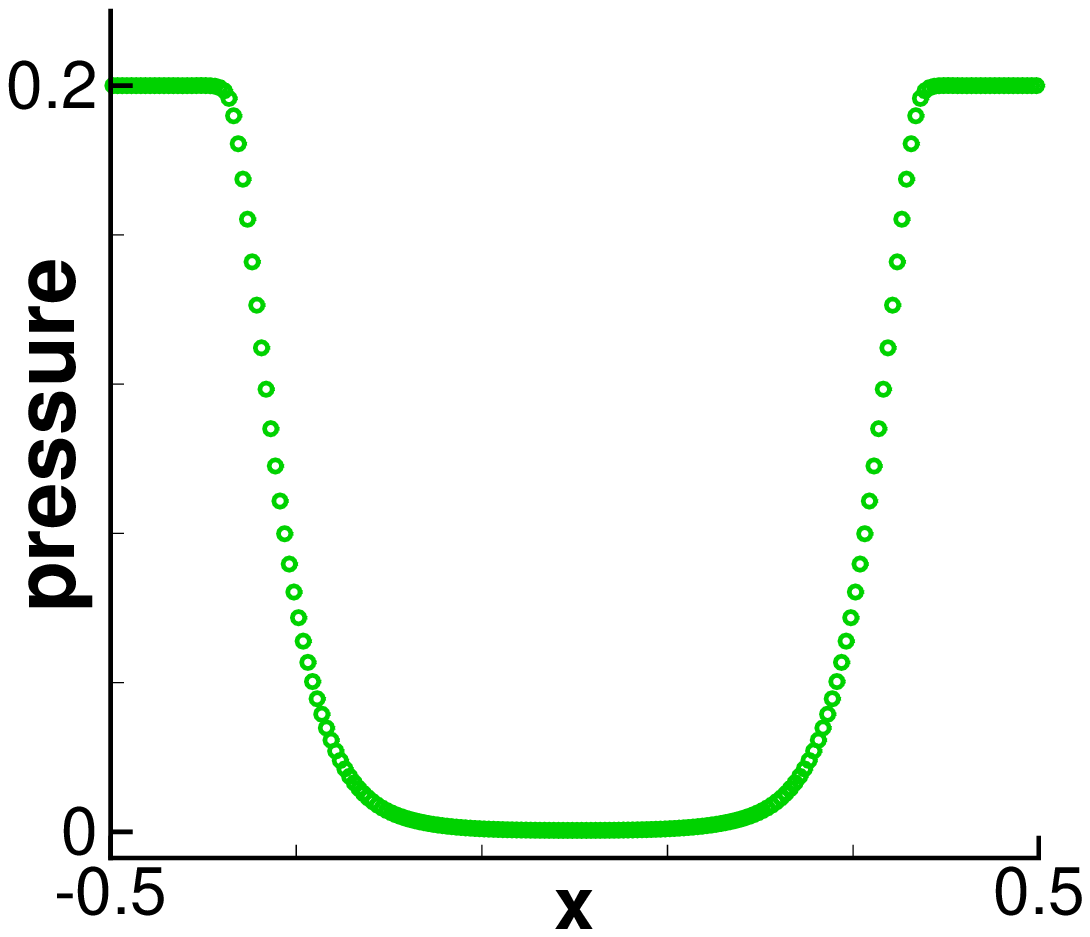}
        \\ {\small Pressure Plot}
    \end{minipage}
    
    \vspace{0.5em}
    
    \begin{minipage}[b]{0.45\textwidth}
        \centering
        \includegraphics[width=\textwidth]{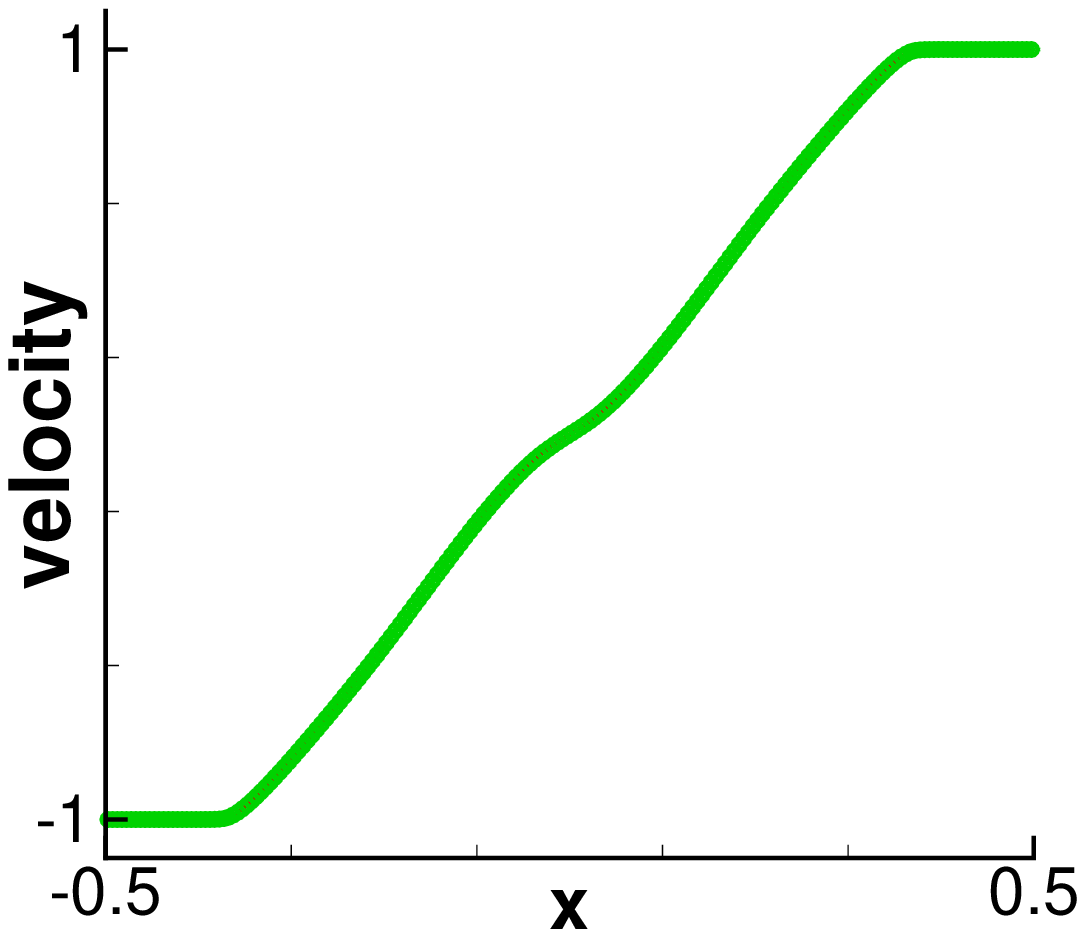}
        \\ {\small Velocity Plot}
    \end{minipage}

    \caption{Example \ref{ex: 1D_rarefaction}: Approximate solution of the 1D double rarefaction problem using 
fifth order WENO with the sweeping limiter.  Green symbols are approximate results.}
    \label{fig:1D_rarefaction}
\end{figure}

\end{example}

\begin{example}[1D Sedov blast]
\label{ex: 1D_Sedov}
We consider the one-dimensional Sedov blast wave of (\ref{eq: 1D_Euler}), which results in low density and pressure. The initial condition resembles a delta function, where density is 1, velocity is 0, and total energy is $10^{-12}$ everywhere except the center cell containing energy constant 3,200,000/$\Delta x$.  The boundary conditions are outflow. We run our numerical test on $x = [-2,2]$ with 800 points until final time $T = 0.001$. We use CFL = 1.2 to produce negative pressure values. Even with a larger CFL, our results in Figure \ref{fig:1D_sedov} closely match the exact solution.

\begin{figure}[h!]
    \centering
    \begin{minipage}[b]{0.45\textwidth}
        \centering
        \includegraphics[width=\textwidth]{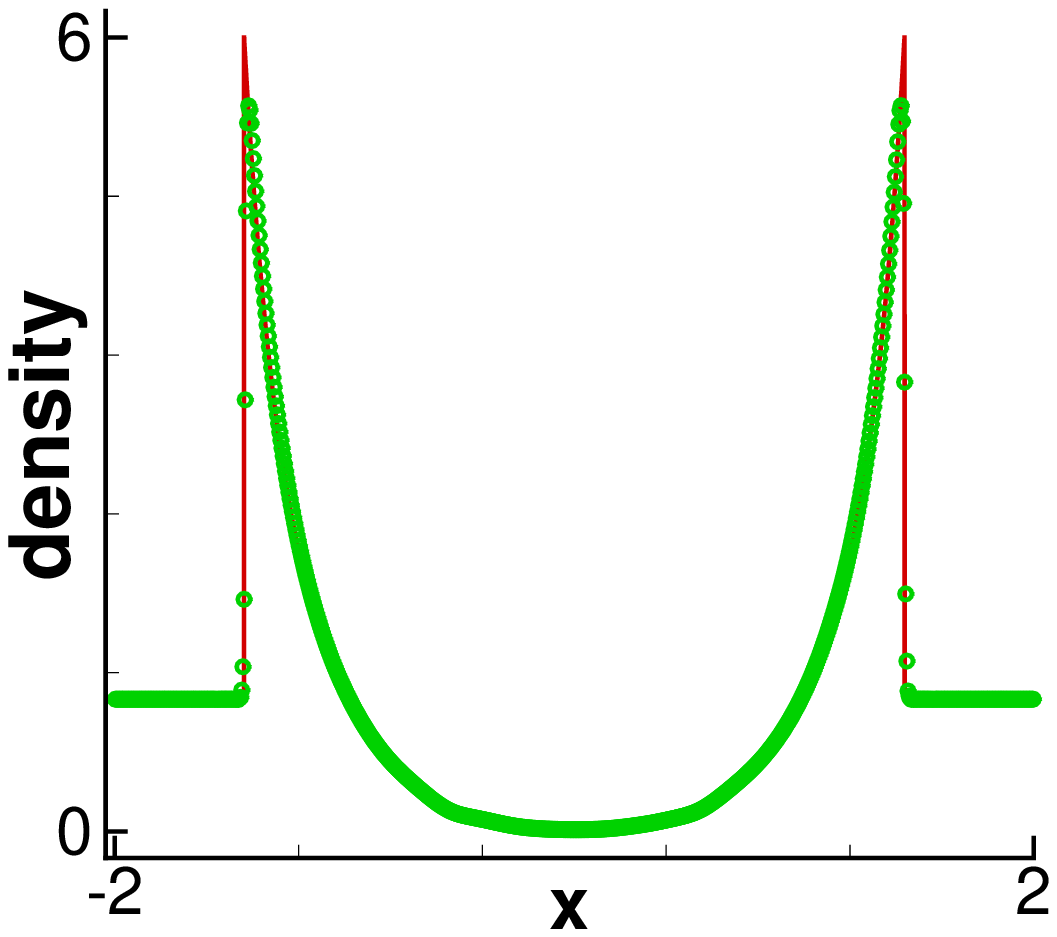}
        \\ {\small Density Plot}
    \end{minipage}
    \hfill
    \begin{minipage}[b]{0.45\textwidth}
        \centering
        \includegraphics[width=\textwidth]{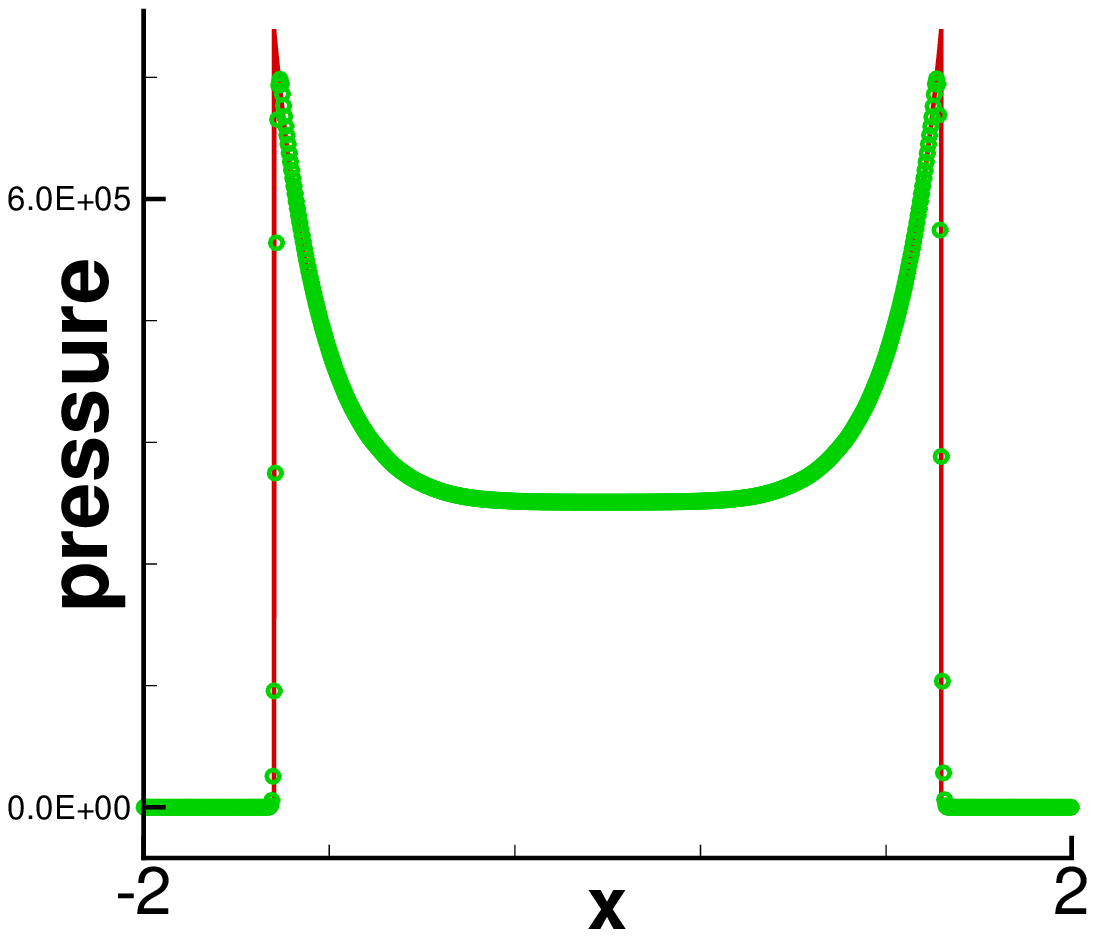}
        \\ {\small Pressure Plot}
    \end{minipage}

    \vspace{0.5em}

    \begin{minipage}[b]{0.45\textwidth}
        \centering
        \includegraphics[width=\textwidth]{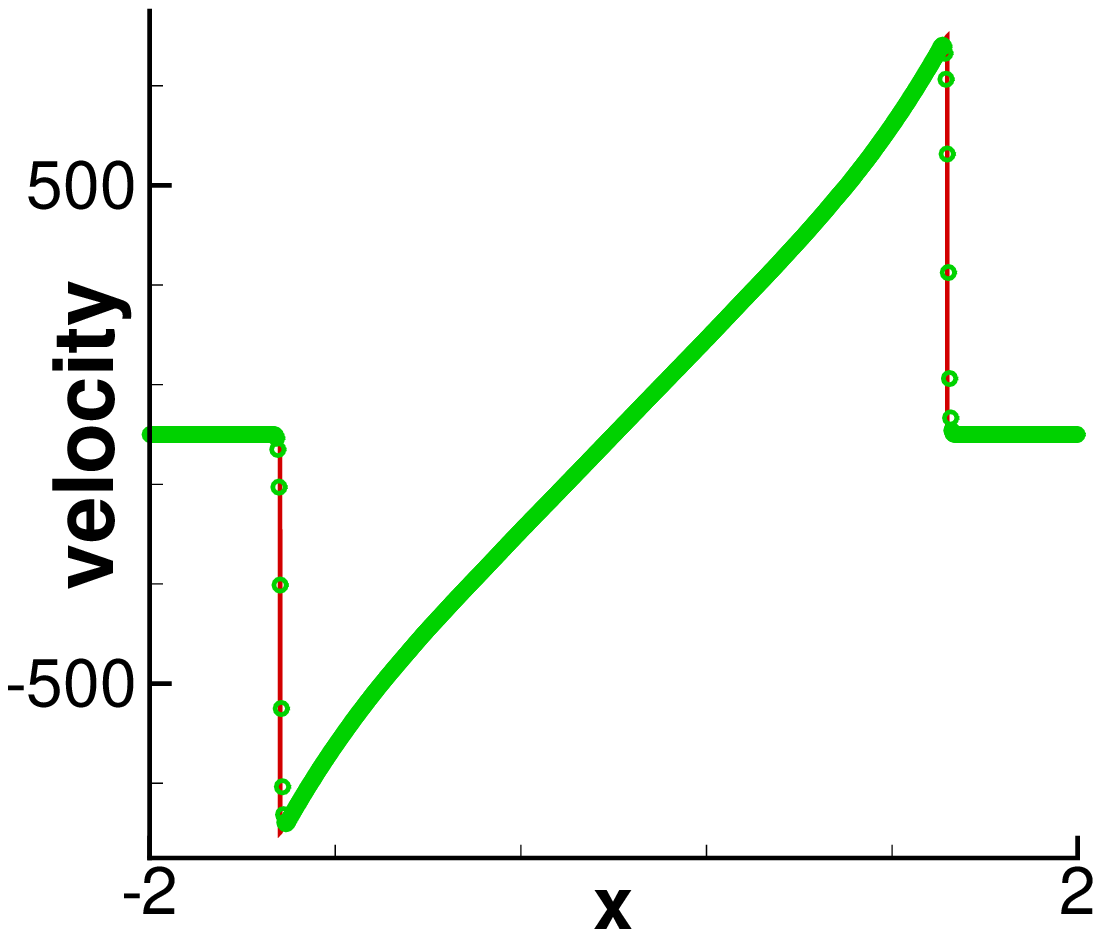}
        \\ {\small Velocity Plot}
    \end{minipage}

    \caption{Example \ref{ex: 1D_Sedov}: Numerical results  1D Sedov blast test case using fifth order WENO with sweeping limiter. Green symbols are approximate results. Red line is exact solution. }
    \label{fig:1D_sedov}
\end{figure}
    
\end{example}

\begin{example}[2D vortex evolution accuracy]
\label{ex: 2D_vortex}
For an accuracy test, we solve the 2D Euler equations (\ref{eq: 2D_Euler}) while considering an isentropic vortex perturbation in $(u,v)$ centered at \((x_0, y_0) = (5,5)\) on $\Omega = [0,10] \times [0,10]$. The mean flow $\rho = u = v = p = 1$ is summed with a vortex pertubation given by 
\[
(\delta u,\delta v) = - \frac{\Gamma}{2\pi} 
\exp\!\left(\tfrac{1-r^2}{2}\right)\,(-\bar{y}, \bar{x}),
\]
\[\delta T = \frac{(\gamma - 1)\Gamma^2}{8\pi^2\gamma} \exp(1 - r^2),\]
and
\[
r^2 = \bar{x}^2 + \bar{y}^2\]
where $\Gamma$ is the vortex strength and $(\bar{x},\bar{y}) = (x - x_0, y - y_0)$. $S = p/\rho^\gamma$ remains unperturbed. This is the same test as \cite{Zhang2012Positivity-preservingEquations}, where $\Gamma = 10.0828$. The exact solution is the advection of the vortex by the mean flow. Using periodic boundaries, we run our numerical test with increasing points, as shown in Table \ref{table: convergence_table} until the final time $T = 0.01$ with CFL = 0.5. We can clearly see the fifth-order accuracy. For a visual representation of the results, we plot a slice through $x = 5$ for the pressure and density with $180 \times 180$ points in space, as shown in Figure \ref{fig:vortex_accuracy}. We also include the sweeping data for this case in Table \ref{table: sweeping_info}.

\begin{figure}[h!]
    \centering
    \begin{minipage}[b]{0.49\textwidth}
        \centering
        \includegraphics[width=\textwidth]{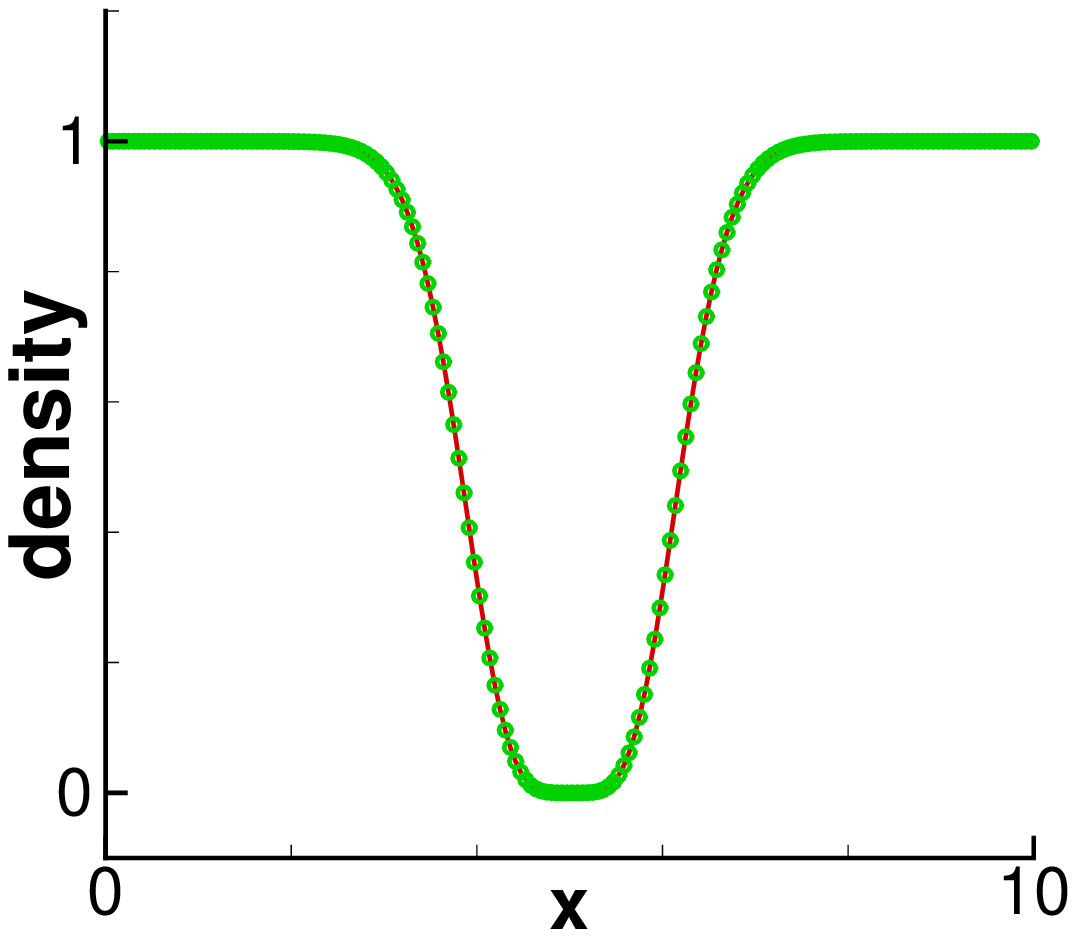}
        \\ {\small Density Plot}
    \end{minipage}
    \hfill
    \begin{minipage}[b]{0.49\textwidth}
        \centering
        \includegraphics[width=\textwidth]{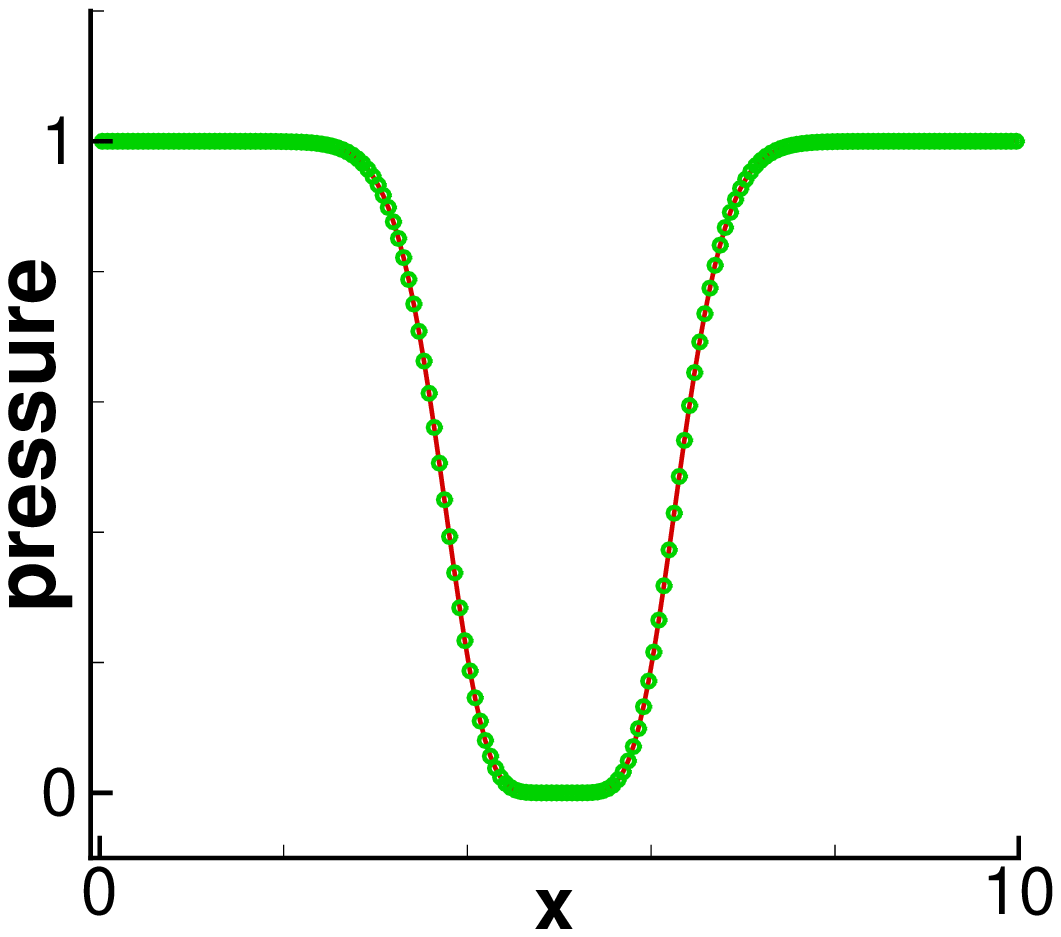}
        \\ {\small Pressure Plot}
    \end{minipage}

     \caption{Example \ref{ex: 2D_vortex}: Numerical results for 2D vortex evolution test case using fifth order WENO with sweeping limiter. Green symbols are approximate results. Red line is exact solution.}
    \label{fig:vortex_accuracy}
\end{figure}

\begin{table}[h!]

\centering
\caption{Convergence Data}
\begin{tabular}{c c c c c}
\hline
Number of points & $L_1$ Error & $L_1$ Order & $L_\infty$ Error & $L_\infty$ Order \\
\hline
$45 \times 45$  & $4.0672 \times 10^{-5}$ & 0.0000 & $2.4255 \times 10^{-3}$ & 0.0000 \\
$90 \times 90$  & $4.6604 \times 10^{-6}$ & 3.1255 & $4.2309 \times 10^{-4}$ & 2.5192 \\
$180 \times 180$& $1.7786 \times 10^{-7}$ & 4.7117 & $2.2275 \times 10^{-5}$ & 4.2475 \\
$360 \times 360$ & $4.1597 \times 10^{-9}$ & 5.4181 & $4.9210 \times 10^{-7}$ & 5.5003 \\
\hline
\label{table: convergence_table}
\end{tabular}
\end{table}

\end{example}

\begin{example}[2D Sedov blast]
\label{ex: 2D_Sedov}
We use the upper-right-quarter computational domain [0,1.3] $\times$ [0,1.3] to solve equations (\ref{eq: 2D_Euler}). The initial condition resembles a delta function, where density is 1, velocity is 0, total energy is $10^{-12}$, and the lower left corner cell contains an energy constant $0.244816/(\Delta x \Delta y)$. The left and bottom boundaries are reflective, and the top and right boundaries are outflow. Using sweeping, we obtain results that are comparable to the exact solution and positivity-preserving method of \cite{Zhang2012Positivity-preservingEquations}. With $640 \times 640$ points in space and CFL = 0.5, we compute the solution at the final time $T = 1$. The results are shown in Figure \ref{fig:2D_sedov}.

\begin{figure}[h!]
    \centering
    \begin{minipage}[b]{0.9\textwidth}
        \centering
        \includegraphics[width=0.8\textwidth]{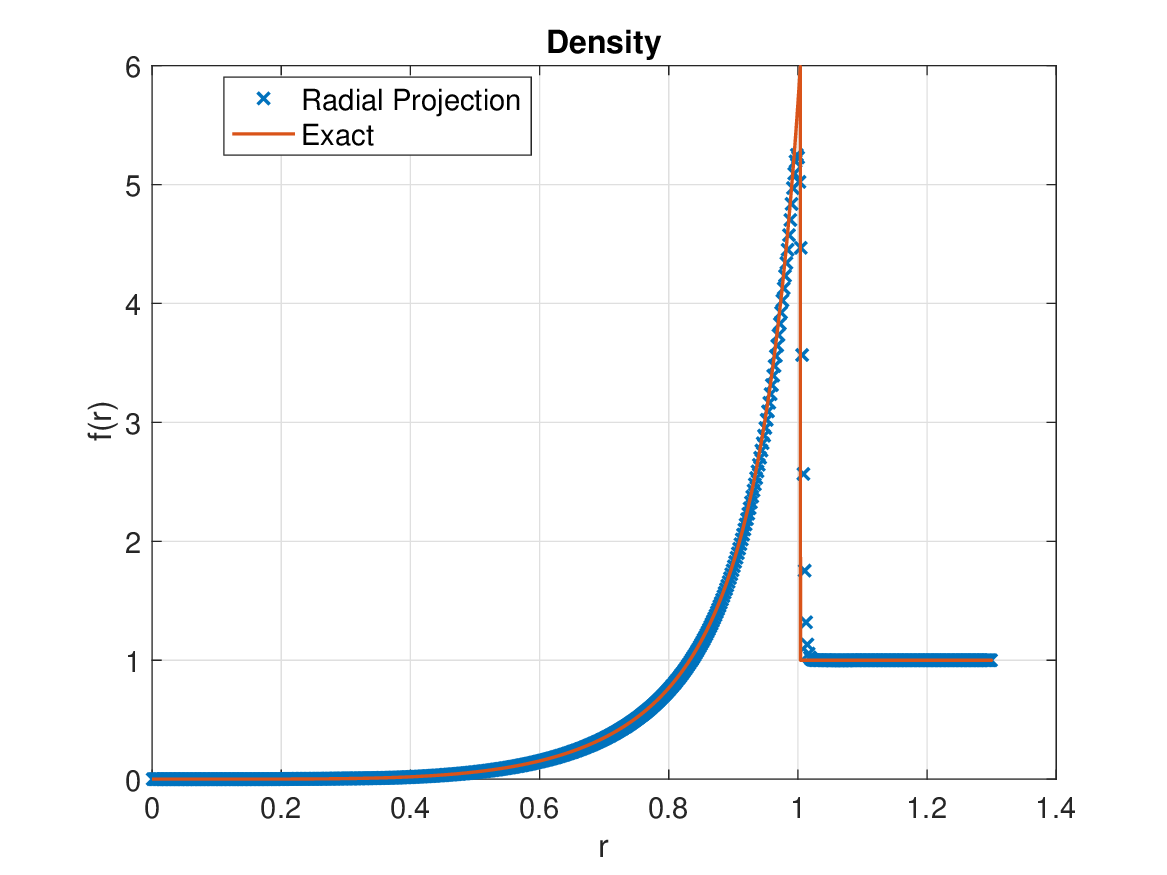}
        \\ {\small Density Radial Projection}
    \end{minipage}

    \vspace{1em}
    \begin{minipage}[b]{0.9\textwidth}
        \centering
        \includegraphics[width=\textwidth]{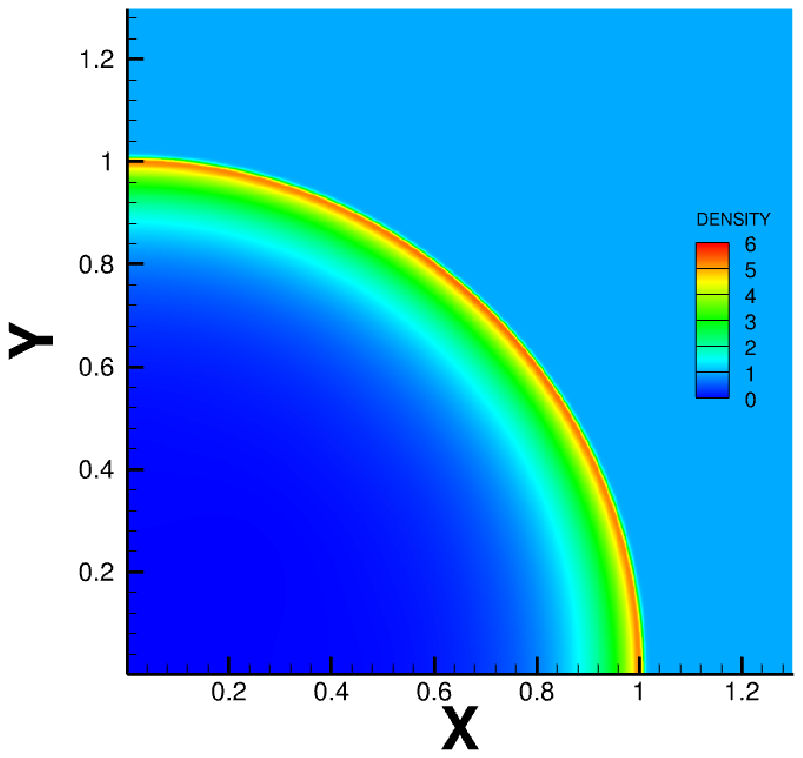}
        \\ {\small Density Contour}
    \end{minipage}

    \caption{Example \ref{ex: 2D_Sedov}: Numerical results for the 2D Sedov blast wave using fifth order WENO with sweeping limiter.  }
    \label{fig:2D_sedov}
\end{figure}

\end{example}

\begin{example}[Mach 2000 robustness test]
\label{ex: Mach_2000}
   To demonstrate the robustness of the sweeping algorithm, we run the same high Mach problem as \cite{Zhang2010OnMeshes} and \cite{Zhang2012Positivity-preservingEquations}. We solve (\ref{eq: 2D_Euler}) with $\gamma = 5/3$. The initial condition is 
   \[(\rho, u, v, p) = (0.5,0,0,0.4127).\]  For the left boundary condition along $x = 0$, we have 

\[\begin{cases}
    (\rho, u, v, p) = (5,800,0,0.4127) & \text{if } y \in [-0.05,0.05],\\
    (\rho, u, v, p) = (5,0,0,0.4127) & \text{if } y \not \in [-0.05,0.05].
\end{cases}\] All other boundaries are outflow. We run our numerical test on $\Omega = [0,1] \times [0,0.25]$ with $800\times 400$ points in space until the final time $T = 0.001$ with CFL = 0.25. Using a jet speed of 800, the high Mach number results in many negative pressure and density values that are corrected by sweeping. The high number of total sweeps is presented in Table \ref{table: sweeping_info}. However, only a maximum of two sweeps were required at any given RK step. The results for density and pressure are shown in Figure \ref{fig:mach_2000}.

\begin{figure}[h!]
    \centering
    \begin{minipage}[b]{0.9\textwidth}
        \centering
        \includegraphics[width=\textwidth]{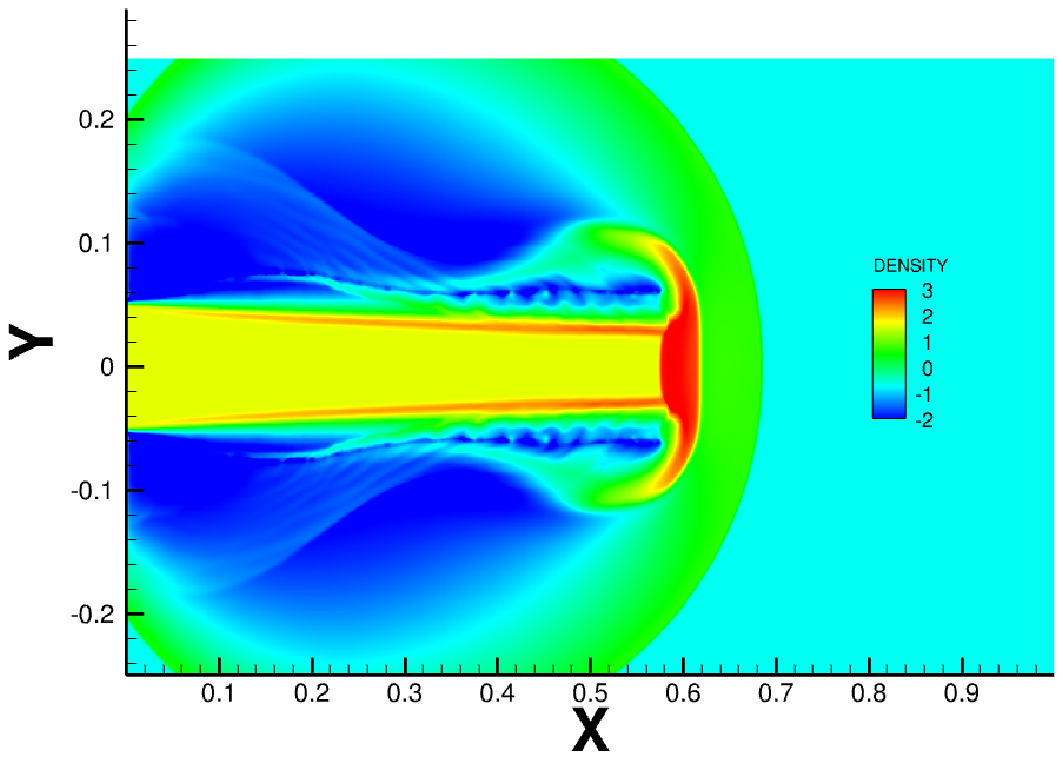}
        \\ {\small Density Plot}
    \end{minipage}

    \vspace{1em}

    \begin{minipage}[b]{0.9\textwidth}
        \centering
        \includegraphics[width=\textwidth]{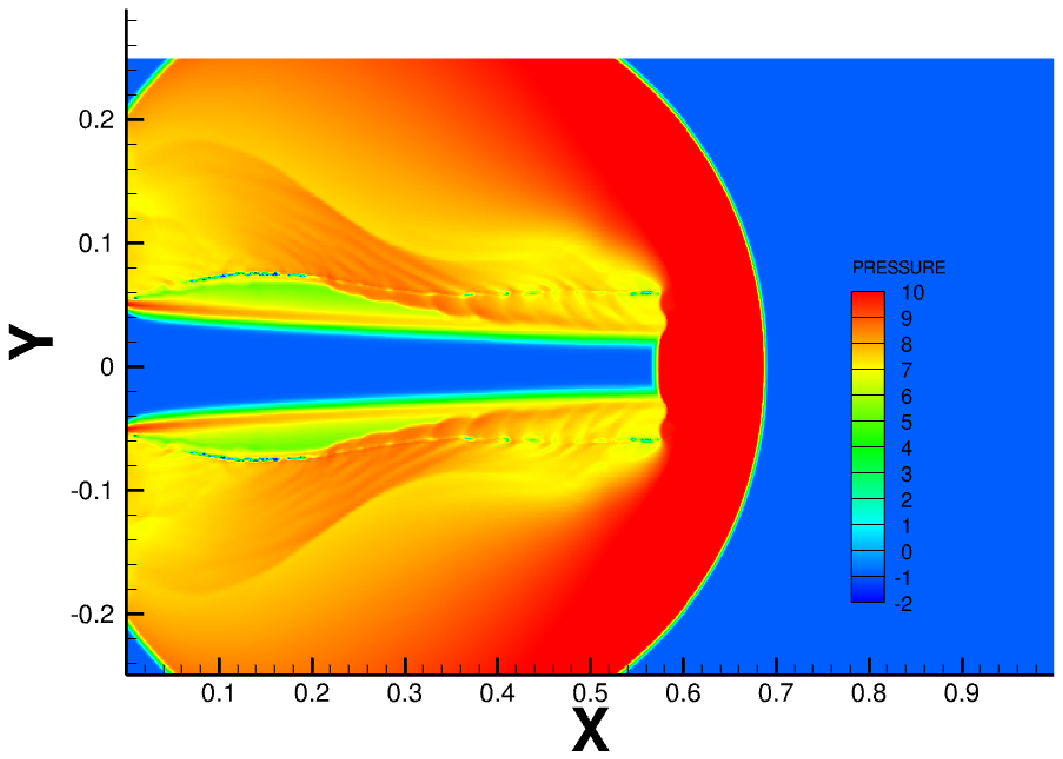}
        \\ {\small Pressure Plot}
    \end{minipage}

     \caption{Example \ref{ex: Mach_2000}: Color contour of results on a logarithmic scale for Mach 2000 test case using 
fifth order WENO with sweeping limiter. }
    \label{fig:mach_2000}
\end{figure}

\end{example}

\begin{example}[Shock diffraction: High CFL]
\label{ex: Shock_Diffraction}
We consider the low-density and low-pressure shock diffraction problem of ($\ref{eq: 2D_Euler})$. The domain is the union of $[0,1] \times [6,11]$ which contains the pure, right-moving Mach 5.09 shock along $x = 0.5$ and $[1,13] \times [0,11]$ which contains ambient gas $(\rho, u, v, p) = (1.4,0,0,1)$. The boundary condition is outflow along the top, right, and bottom of $[1,13] \times [0,11]$. The boundary condition is inflow on the left and reflective on the bottom of $[0,1] \times [6,11]$. Along the wall $\{0\}\times[0,6]$, reflective boundaries are used. We run our numerical test on $\Omega = [0,13] \times [0,11]$ with $1040 \times 880$ points until the final time $T = 2.3$. We use a high CFL = 0.9 to produce negative pressure values. Due to the high CFL, we see small fluctuations in the solution. Nevertheless, the plot shows key features, even at the instability threshold. The results are shown in Figure \ref{fig:shock_diffraction}.

\begin{figure}[H]
    \centering
    \begin{minipage}[b]{0.45\textwidth}
        \centering
        \includegraphics[width=\textwidth]{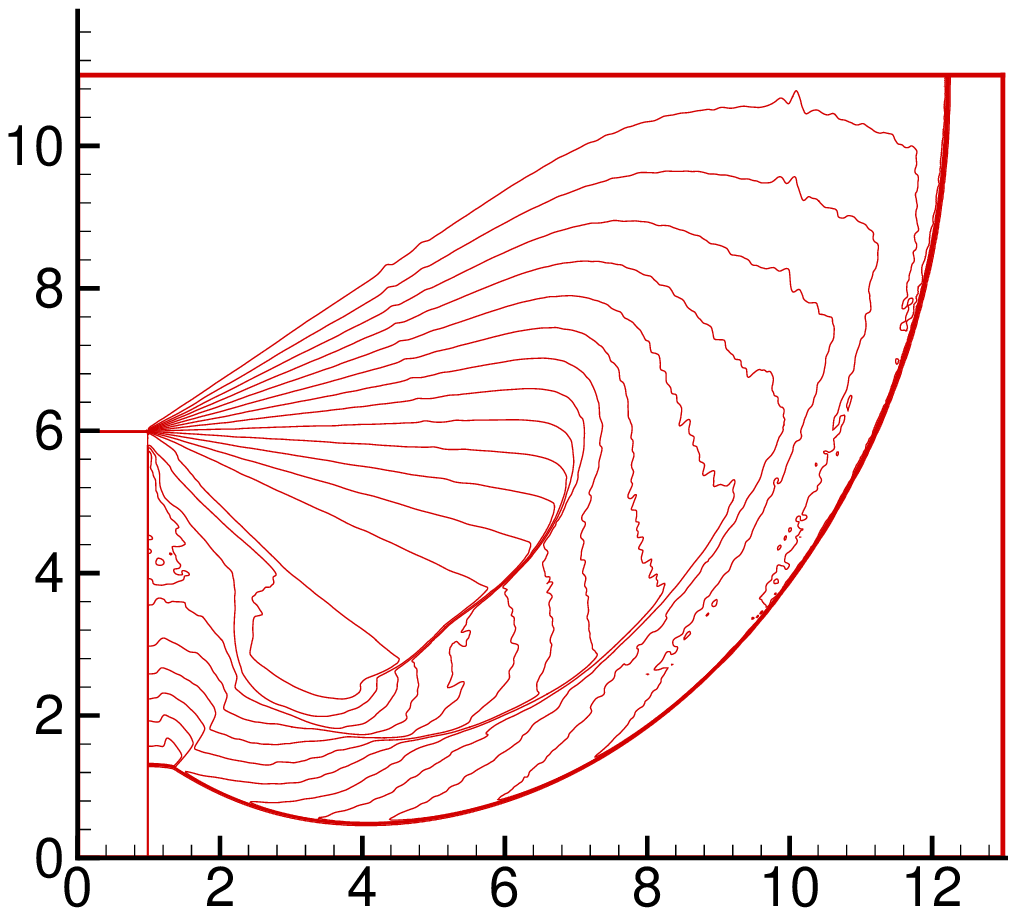}
        \\ {\small Density Plot}
    \end{minipage}
    \hfill
    \begin{minipage}[b]{0.45\textwidth}
        \centering
        \includegraphics[width=\textwidth]{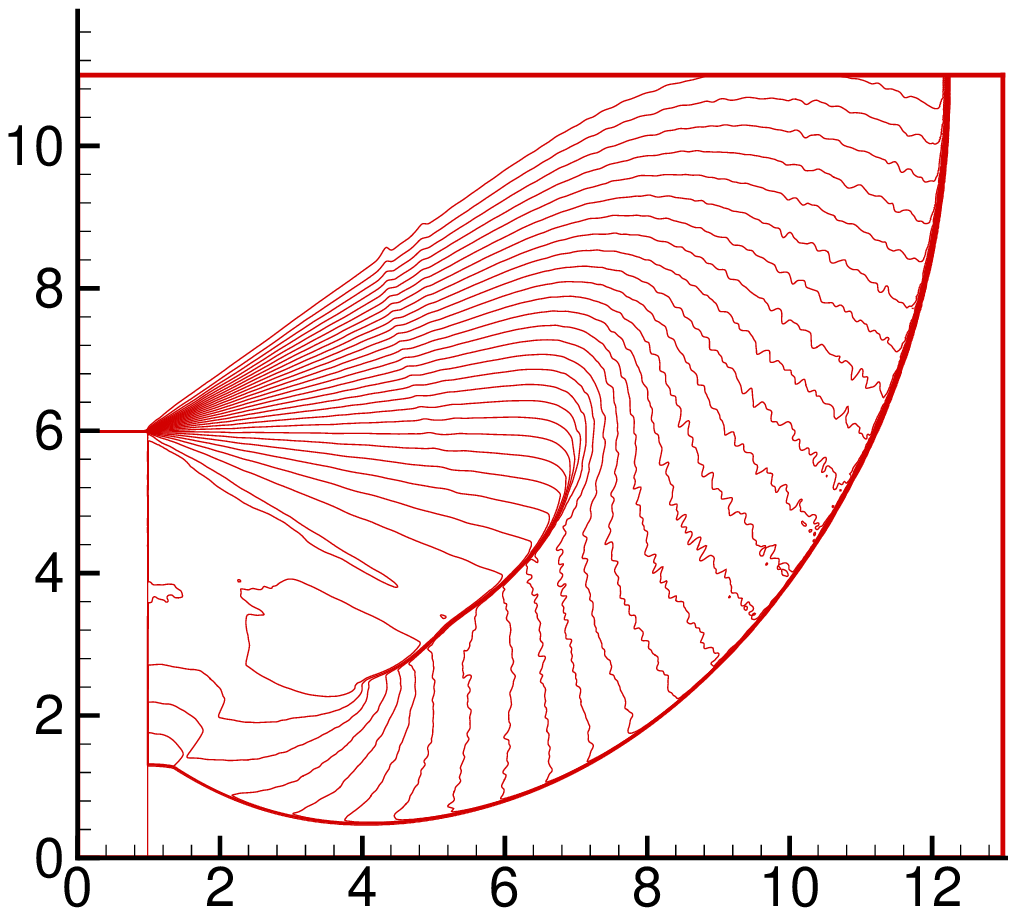}
        \\ {\small Pressure Plot}
    \end{minipage}

    \caption{Example \ref{ex: Shock_Diffraction}: Contour plot of results for shock diffraction test case using fifth order WENO with sweeping limiter. Left plot: 20 equally spaced contour lines between 0.066227 and 7.0668. Right plot: 40 contour lines between 0.091 and 42.}
    \label{fig:shock_diffraction}
\end{figure}
\end{example}
\begin{example}[Reactive Euler detonation: High CFL]
\label{ex: detonation_diffraction}
  The reactive Euler equations are an extension of the standard Euler equations for modeling chemical reactions or detonations. The governing equations are written as 
\begin{equation}
    \mathbf{w}_t + \mathbf{f}(\mathbf{w})_x + \mathbf{g}(\mathbf{w})_y = \mathbf{s}(\mathbf{w}), \; t \geq 0, \quad (x,y) \in \mathbb{R}^2,
\end{equation}

\[
\mathbf{w} = \begin{pmatrix}\rho\\m\\n\\E\\\rho Y\end{pmatrix}, \quad 
\mathbf{f}(\mathbf{w}) = \begin{pmatrix}
    m\\ \rho u^2 + p\\ \rho u v \\ (E + p) u\\ \rho u Y
\end{pmatrix}, \quad 
\mathbf{g}(\mathbf{w}) = \begin{pmatrix} 
    n \\ \rho u v\\ \rho v^2 + p \\ (E + p) v \\ \rho v Y
\end{pmatrix}, \quad 
\mathbf{s}(\mathbf{w}) = \begin{pmatrix}
    0 \\ 0 \\ 0 \\ 0 \\ \omega
\end{pmatrix},
\]

where 
\[
E = \tfrac{1}{2} \rho u^2 + \tfrac{1}{2} \rho v^2 + \tfrac{p}{\gamma - 1} + \rho \mathcal{Q} Y,
\]
\[
m = \rho u, \quad n = \rho v, \quad \omega = -K \rho Y \exp\!\left(-{E_a}/{\tilde{T}}\right), \tilde{T} = p/\rho.
\]
Here, $\omega$ is in the Arrhenius kinetic form, $\gamma$ is the heat ratio, $\tilde{T}$ is the temperature, $Y$ is the reactant mass fraction, and $\mathcal{Q}$ is the heat release of the reaction. Solving the energy equation for pressure, we obtain
\[
p = (\gamma - 1)\left(\frac{m^2 + n^2}{2\rho} - \mathcal{Q}(\rho Y)\right).
\]
The sum of concave functions is concave. For $\rho > 0$, $\frac{m^2 + n^2}{2\rho}$ is clearly concave, and $\mathcal{Q}(\rho Y)$ is linear for a constant $\mathcal{Q}$. Thus, sweeping is applied to the modified pressure function. The domain is  $[0,1] \times [2,5] \cup [1,5] \times [0,5]$. The initial condition is given by 

\[\begin{cases}
    (\rho, u ,v, E, Y) = (11,6.18,0,970,1), &\text{if } x < 0.5,\\
    (\rho, u ,v, E, Y) = (1,0,0,55,1), &\text{if } x \geq 0.5
\end{cases}\]
At $x = 0,$ we have the boundary condition 
\[(\rho, u ,v, E, Y) = (11,6.18,0,970,1),\]
while the remaining boundaries are reflective. We run our numerical test on $\Omega = [0,5] \times [0,5]$ with $400 \times 400$ points until the final time $T = 0.6$. We use CFL = 0.89 to produce negative pressure values. We use constants $K = 2566.4$, $\mathcal{Q} = 50$, $E_a = 50$, and $\gamma = 1.2$. Even with a high CFL number at the instability threshold of the underlying scheme, the sweeping method prevents blow-ups of the numerical solution and preserves key features. The results are shown in Figure \ref{fig:reactive_euler}.

\begin{figure}[H]
    \centering
    \begin{minipage}[b]{0.45\textwidth}
        \centering
        \includegraphics[width=\textwidth]{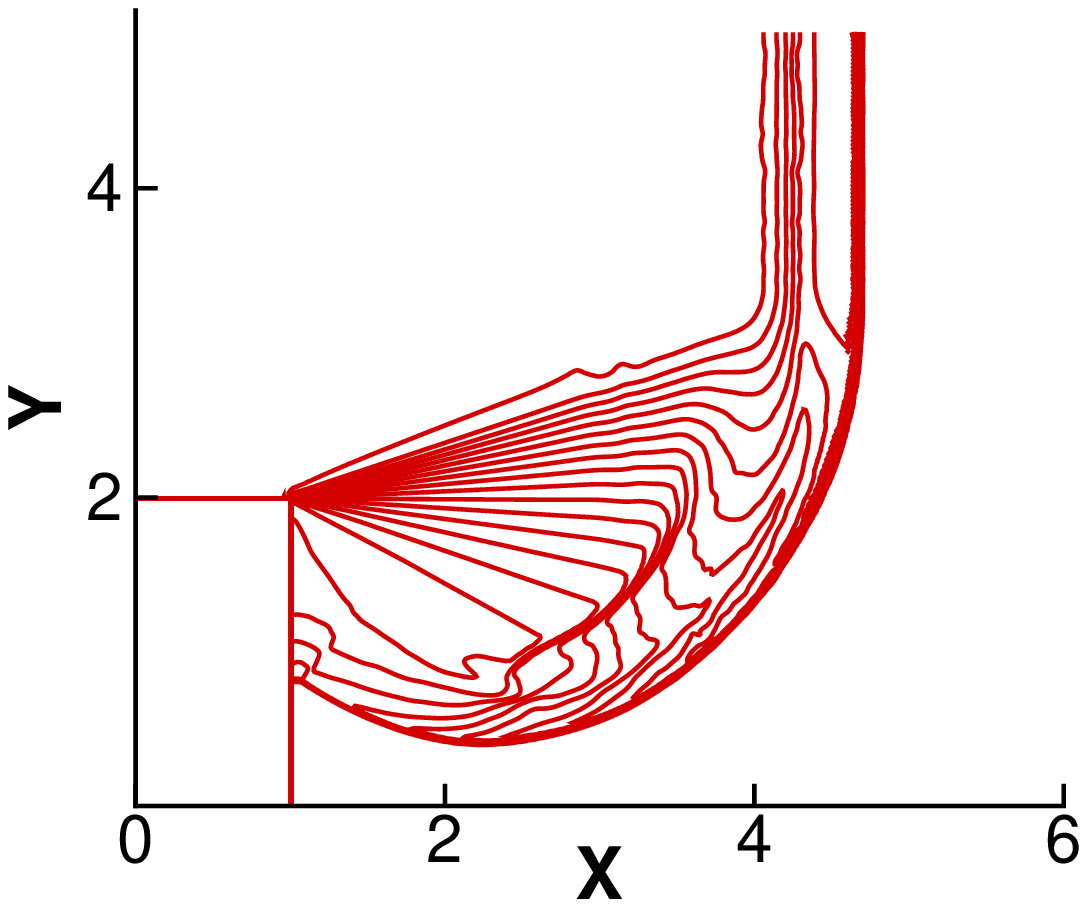}
        \\ {\small Density Contour Lines}
    \end{minipage}
    \hfill
    \begin{minipage}[b]{0.45\textwidth}
        \centering
        \includegraphics[width=\textwidth]{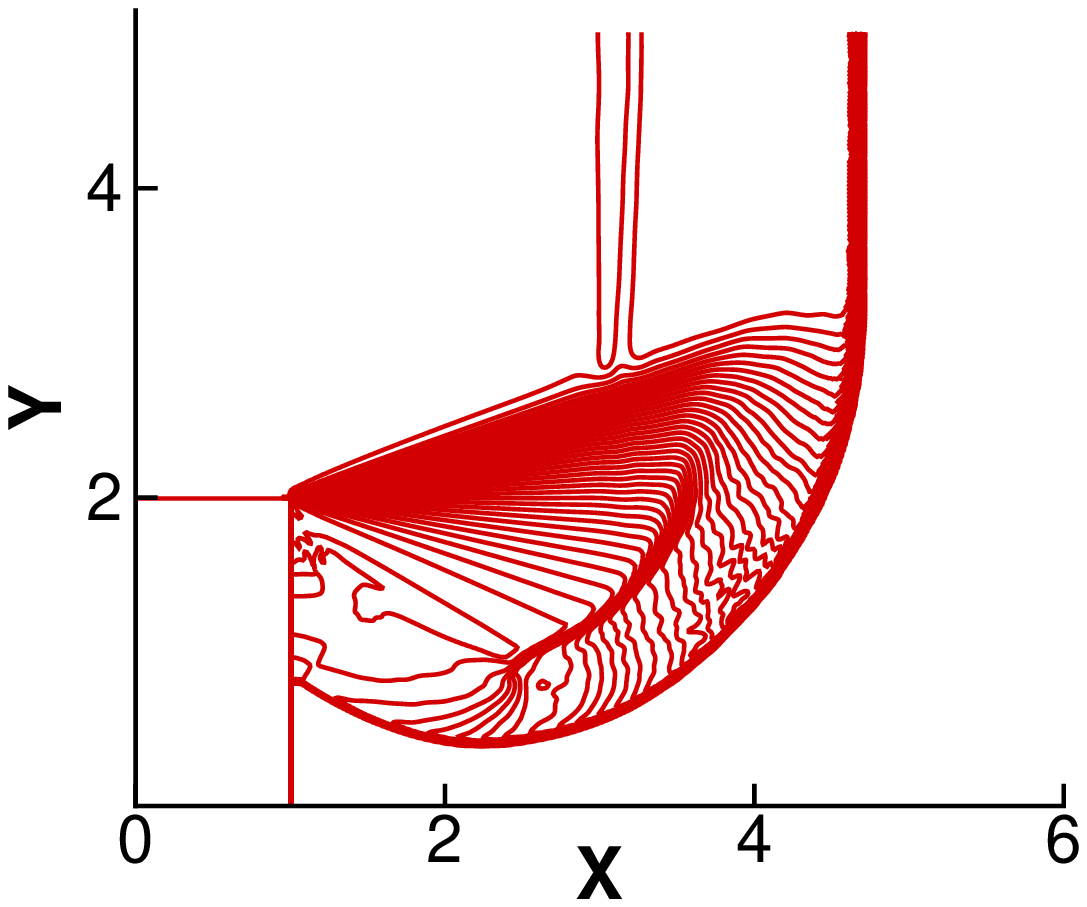}
        \\ {\small Pressure Contour Lines}
    \end{minipage}

    \vspace{0.5em}

    \begin{minipage}[b]{0.45\textwidth}
        \centering
        \includegraphics[width=\textwidth]{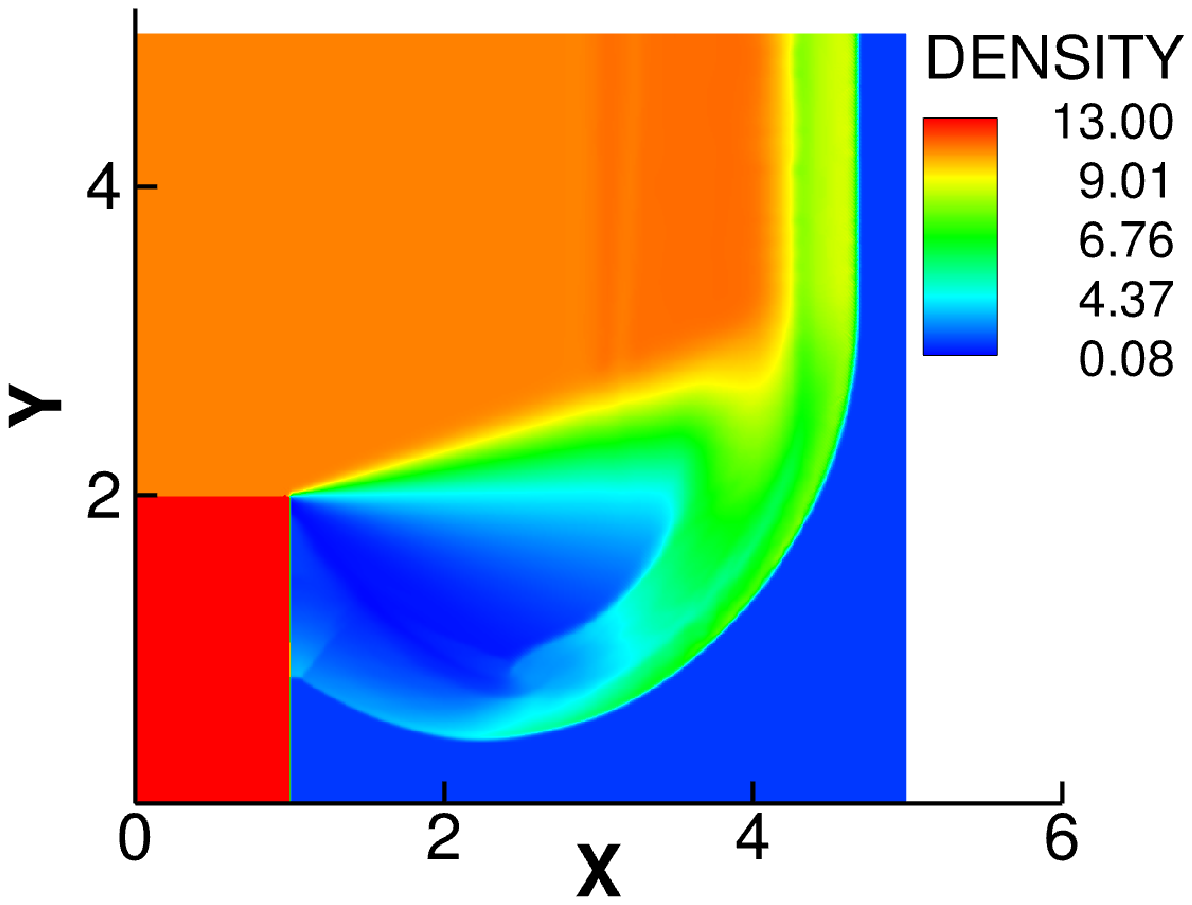}
        \\ {\small Density Color Contour}
    \end{minipage}
    \hfill
    \begin{minipage}[b]{0.45\textwidth}
        \centering
        \includegraphics[width=\textwidth]{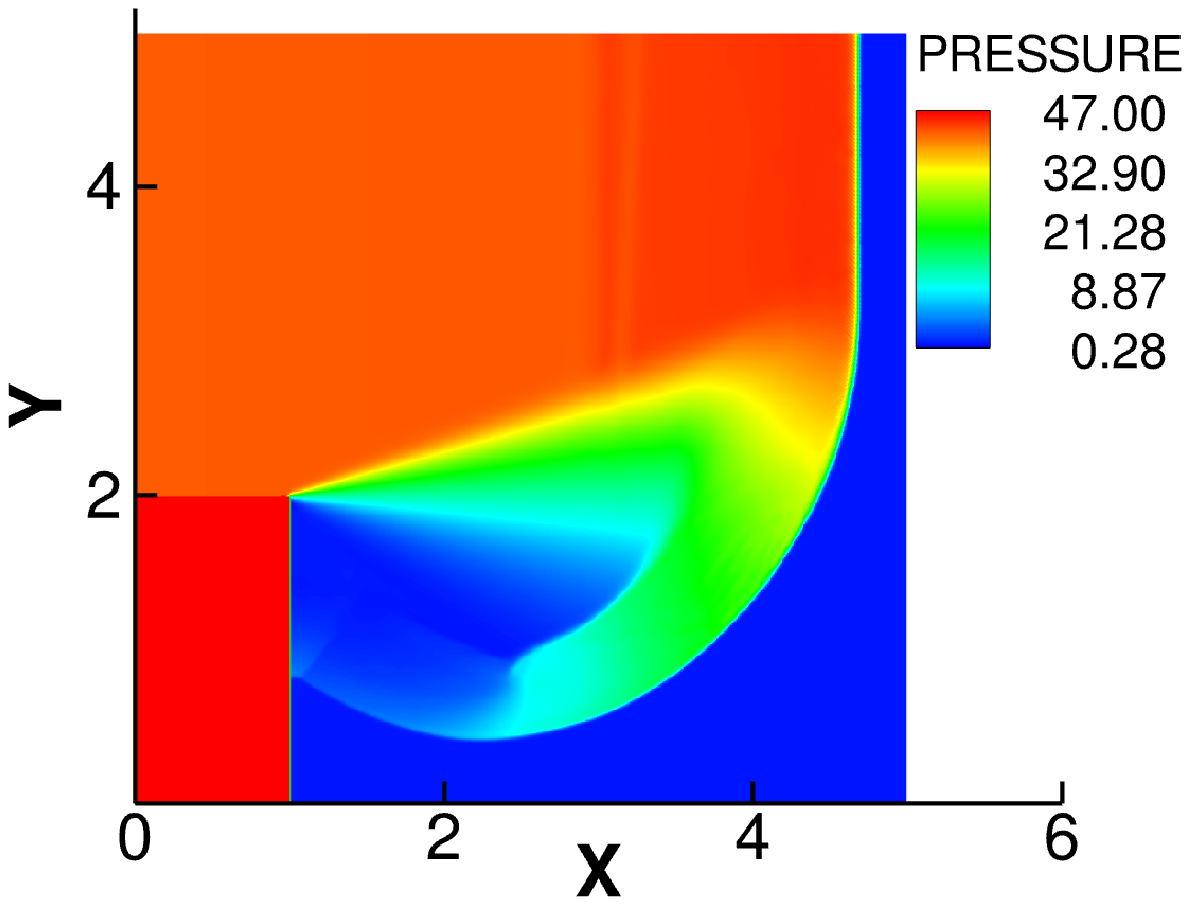}
        \\ {\small Pressure Color Contour}
    \end{minipage}

    \caption{Example \ref{ex: detonation_diffraction}: Approximate solution of detonation diffraction test case of reactive Euler equations using fifth order WENO with sweeping limiter. Upper left plot: 20 contour lines between 0 and 13. Upper right plot: 40 contour lines between 0 and 43.}
    \label{fig:reactive_euler}
\end{figure}
\end{example}
Here, we provide a table of the maximum number of sweeps required at any given RK step. We also provide the average and total number of sweeps required by the final time. A full sweep includes the forward and backward direction. In most cases, only one forward and backward sweep is necessary for positivity preservation. We reiterate that the density sweep requires at most one forward and backward sweep, thus we only provide the total number of density sweeps until final time.

\begin{table}[h!]

\centering
\caption{Pressure Sweeping Data}
\begin{tabular}{lcccc}
\hline
\label{table: sweeping_info}
 & 1D Rarefaction & 1D Sedov & Vortex Accuracy & 2D Sedov \\
\hline
Max Sweep     & 1 & 1 & 1 & 1 \\
Average Sweep & 1 & 1 & 1 & 1 \\
Total Sweep   & 2 & 1 & 8 & 3 \\
\hline
\end{tabular}

\begin{tabular}{lccc}
\hline
 & Mach 2000 & Shock Diffraction & Reactive Euler \\
\hline
Max Sweep     & 2 & 2 & 1 \\
Average Sweep & 1.9106 & 1.0004 & 1 \\
Total Sweep   &  18185 & 2376 & 286 \\
\hline
\end{tabular}
\end{table}

\begin{table}[h!]
\centering
\caption{Density Sweeping Data}
\begin{tabular}{lcccc}
\hline
\label{table: den_sweeping_info}
 & 1D Rarefaction & 1D Sedov & Vortex Accuracy & 2D Sedov \\
\hline
Total Sweep   & 0 & 1 & 0 & 0  \\

\hline
\end{tabular}

\begin{tabular}{lccc}
\hline
 & Mach 2000 & Shock Diffraction & Reactive Euler \\
\hline
Total Sweep   &  0 & 890 & 22 \\
\hline
\end{tabular}
\end{table}
\newpage
\section{Concluding Remarks}
In this work, we provide a conservative, positivity-preserving sweeping procedure for the pressure function in the Euler equations of gas dynamics as a post-processing step for both finite difference and finite volume methods. We show that the procedure is conservative and positivity-preserving. Our results for fifth order finite difference WENO demonstrate that the sweeping procedure is high-order and robust with only one or two sweeps needed at any given Runge-Kutta time-stepping stage,
and high order accuracy is always maintained. 

We aim to extend the sweeping procedure to the Navier-Stokes equations and other applications where
positivity for concave functions of the conserved variables is required. Other applications include 
magnetohydrodynamics using implicit algorithms. 
\label{section: concluding_remarks}
\section{Acknowledgements}
The research of the first author is partially supported by the NSF
Graduate Research Fellowship Program. The research of the second author is partially supported
by NSF grant DMS-2309249. 

\bibliographystyle{plainnat}

\end{document}